\documentclass[graybox,envcountsame,envcountsect]{SNmult}
\usepackage{type1cm}  %
\makeatletter\let\@font@warning\@gobble\makeatother
\usepackage{tikz}
\usetikzlibrary{decorations.markings, arrows.meta}
\usepackage{comment}

\usepackage{amsmath,amsthm,amsxtra,amsbsy,accents}
\usepackage{latexsym,amssymb}
\usepackage{mathtools}
\usepackage[scr=boondoxo]{mathalpha}
\usepackage{bm}
\usepackage[normalem]{ulem}  %
\usepackage{lmodern,microtype}

\usepackage{enumitem}

\usepackage{upgreek}

\usepackage[only,llbracket,rrbracket]{stmaryrd}
\SetSymbolFont{stmry}{bold}{U}{stmry}{m}{n}

\DeclareMathAlphabet{\mathup}{OT1}{\familydefault}{m}{n}
\newcommand{\dx}[1]{\mathop{}\!\mathup{d} #1}
\newcommand{\pderiv}[3][]{\frac{\mathop{}\!\mathup{d}^{#1} #2}{\mathop{}\!\mathup{d} #3^{#1}}}

\newcommand{\eps}{\varepsilon}
\DeclarePairedDelimiter{\abs}{\lvert}{\rvert}
\DeclarePairedDelimiter{\norm}{\lVert}{\rVert}
\DeclarePairedDelimiter{\bra}{(}{)}
\DeclarePairedDelimiter{\pra}{[}{]}
\DeclarePairedDelimiter{\set}{\{}{\}}
\DeclarePairedDelimiter{\skp}{\langle}{\rangle}
\usepackage{ifthen}
\newlength{\leftstackrelawd}
\newlength{\leftstackrelbwd}
\def\leftstackrel#1#2{\settowidth{\leftstackrelawd}%
{${{}^{#1}}$}\settowidth{\leftstackrelbwd}{$#2$}%
\addtolength{\leftstackrelawd}{-\leftstackrelbwd}%
\leavevmode\ifthenelse{\lengthtest{\leftstackrelawd>0pt}}%
{\kern-.5\leftstackrelawd}{}\mathrel{\mathop{#2}\limits^{#1}}}
 \def\calB{{\mathcal B}} 
\def\calD{{\mathcal D}} \def\calE{{\mathcal E}} 
 \def\calH{{\mathcal H}} \def\calI{{\mathcal I}}
  \def\calL{{\mathcal L}}
\def\calM{{\mathcal M}}  
\def\calP{{\mathcal P}}  \def\calR{{\mathcal R}}
\def\calS{{\mathcal S}}  \def\calU{{\mathcal U}}
\def\calV{{\mathcal V}}

\def\rmj{{\mathrm j}}

\def\rmD{{\mathrm D}}

 \def\sfB{{\mathsf B}} \def\sfC{{\mathsf C}}
 \def\sfE{{\mathsf E}}

\def\sfM{{\mathsf M}}  
  \def\sfR{{\mathsf R}}

  \def\scrL{{\mathscr  L}}

\newcommand{\R}{\mathbb{R}}

\newcommand{\CE}{{\mathsf{CE}}}
\newcommand{\lCE}{{\overline{\mathsf{CE}}}}

\renewcommand{\div}{\mathop{\mathrm{div}}\nolimits}

\newcommand{\ub}{{\overline u}}
\newcommand{\jb}{\bar \jmath}

\DeclareMathOperator{\tr}{tr}
\DeclareMathOperator{\arsinh}{arsinh}
\DeclareMathOperator{\DtN}{D_{to}N}

\makeatletter
\newcommand{\bnabla}{{\mathpalette\b@nabla\relax}}
\newcommand\b@nabla[2]{%
        \setbox\z@=\hbox{$\m@th#1\bigtriangledown$}%
        \ht\z@.7\ht\z@
        \raise\dp\z@\box\z@
}
\makeatother

\makeatletter
\newcommand{\oset}[3][0ex]{%
  \mathrel{\mathop{#3}\limits^{
    \vbox to#1{\kern-2\ex@
    \hbox{$\scriptstyle#2$}\vss}}}}
\makeatother

\DeclareMathOperator*{\esssup}{ess\,sup}

\newcommand{\FI}{\mathit{F\mkern-4muI}}

\usepackage{trimclip}
\makeatletter
\NewDocumentCommand{\rmjmath}{}{\mathbin{\mathpalette\eplus@\relax\mspace{1mu}}}
\newcommand{\eplus@}[2]{\clipbox{-.5 -.5 0 {.35\height}}{$\m@th#1\rmj$}}
\makeatother

\usepackage[hyphens]{url}

\spnewtheorem{result}{Result}{\bfseries}{\itshape}

\makeatletter
\newcommand{\anabla}{{\mathpalette\a@nabla\relax}}
\newcommand\a@nabla[2]{%
	\setbox\z@=\hbox{$\m@th#1\bigtriangledown$}%
	\ht\z@.7\ht\z@
	\raise\dp\z@\box\z@
}
\makeatother

\usepackage[hidelinks,bookmarks=false,hypertexnames=false]{hyperref}
\hypersetup{
  pdfauthor={Georg Heinze, Mark Peletier, Andr\'e Schlichting},
  pdftitle={Entropic gradient flow structure of the Dirichlet-to-Neumann evolution
        via multidimensional boundary cell problems},
}
\makeatletter
\providecommand*{\toclevel@title}{0}
\providecommand*{\toclevel@author}{0}
\makeatother

\newcommand{\Acc}{\mathcal{A}}

\begin{document}

\motto{\emph{In honour of Giuseppe Savar\'e's long-standing commitment
both to great science and to great colleague-ship}}

\title*{Entropic gradient flow structure of the Dirichlet-to-Neumann evolution}
\subtitle{via multidimensional boundary cell problems}
\titlerunning{Entropic gradient flow structure of the Dirichlet-to-Neumann evolution}

\author{Georg Heinze \and Mark Peletier \and Andr\'e Schlichting}
\authorrunning{Heinze, Peletier, Schlichting}

\institute{Georg Heinze \at Research Group \emph{Partial Differential Equations}, Weierstrass Institute, Mohrenstrasse~39, 10117 Berlin, Germany.
\email{georg.heinze@wias-berlin.de}
\and Mark Peletier \at Department of Mathematics and Computer Science\\
        Eindhoven University of Technology 
 PO Box 513, 5600 MB Eindhoven, The Netherlands
 \email{m.a.peletier@tue.nl}
\and Andr\'e Schlichting \at Institute of Applied Analysis, Ulm University, Helmholtzstrasse~18, 89081 Ulm, Germany.
\email{andre.schlichting@uni-ulm.de}}

\maketitle

\abstract{We derive the Dirichlet-to-Neumann ($\DtN$) evolution on the boundary
$\Gamma=\partial\Omega$ of a bounded smooth domain $\Omega\subset\R^d$ as a
multiscale limit of diffusion equations on $\Omega\cup\Gamma$ with accelerated
interior diffusion. The prelimit systems are gradient flows of the relative
entropy with respect to a combination of Otto- and $\cosh$-type dissipation
potentials. The limiting gradient structure is governed by a
\emph{multidimensional cell problem} relating boundary concentrations to
boundary fluxes, which generalizes a cell problem from diffusion-to-reaction
limits. We show that this cell formula admits an $R$-$R^*$ splitting, exhibiting
the $\DtN$ evolution as a generalized gradient flow driven by the boundary
entropy. Throughout we work with densities that are bounded away from
zero and infinity; this keeps every function space classical and every argument
short, and in this regime all structural features of the limit are
already visible. 
\keywords{Dirichlet-to-Neumann operator $\cdot$ gradient flow $\cdot$ cell problem $\cdot$ multiscale limit $\cdot$ non-equilibrium steady states}}

\newpage

\section{Introduction}

In this paper we bring together three independent strands of research to make a
small but interesting observation. The first strand is Giuseppe's extensive line
of work on gradient flows in all their forms. In particular, in the paper
\emph{Passing to the limit in a Wasserstein gradient flow: from diffusion to
reaction}~\cite{ArnrichMielkePeletierSavareVeneroni2012}, Giuseppe and his
co-authors established the convergence of a diffusion equation to an equation
describing a chemical reaction. At the core of this result is a \emph{cell
problem} which characterizes the limiting relationship between the
concentrations at the two ends of a line segment, on the one hand, and the flux
between these boundary points on the other hand. Later Liero, Mielke, Peletier and Renger
showed~\cite{LieroMielkePeletierRenger2017} that this cell problem admits an
\emph{$R$-$R^*$ structure} with cosh-type dissipation potential $R^*$, implying
that the limiting system describing the chemical reaction is a generalized
gradient flow.

\medskip

The second strand is the fact that on a domain $\Omega\subset \R^d$ with
finite-area boundary the \emph{Dirichlet-to-Neumann} ($\DtN$) operator generates
a Markovian semigroup in the space of non-negative measures on the boundary
$\partial \Omega$. This was independently established by Hsu~\cite{Hsu86} and by Arendt and Ter Elst~\cite{ArendtElst2011}; see also~\cite{SatoUeno65,BezneaVladoiu17,ArendtElst2020,ElstOuhabaz2019}. In particular, Hsu showed that $\DtN$ is the generator of the \emph{boundary process}, a jump process on $\partial \Omega$ defined in terms of diffusion excursions from boundary point to boundary point. In this setup the time is rescaled by the boundary local time, implying that the natural time is no longer visible, rescaled time is random, and the excursions are given zero duration. 

In this paper we also connect the $\DtN$ evolution to classical diffusion, but in a different way: we show that it is a limit of
$\eps$-dependent diffusion evolution equations on $\Omega \cup \partial \Omega$.
For $\eps>0$, these processes have a gradient structure, and the convergence
provides a limiting variational structure that involves a cell problem. This cell
problem turns out to be a direct generalization of the line-segment cell problem 
in~\cite{ArnrichMielkePeletierSavareVeneroni2012}  to the
multidimensional domain~$\Omega$ of this paper.

We further show that this cell problem admits an $R$-$R^*$ splitting, implying that the $\DtN$ evolution is a generalized gradient flow.

\medskip

The third strand of research is Mielke's recent development~\cite{Mielke2023} of
the concept of \emph{Non-Equilibrium Steady States} (NESS) of gradient systems,
and in particular their role in multiscale limits of gradient systems. 
We discuss the derived cell problem and its $R$-$R^*$ structure in the context of this NESS theory in Section~\ref{sec:discussion}. 
This not only motivates the result from a slightly different point of view, but also highlights the role of the Fisher information
as a generator of the relevant function classes,
and thus connects back to Giuseppe's
long-standing interest in Fisher informations and their
properties~\cite{GianazzaSavareToscani2008}.

\medskip
\emph{A simplifying assumption.} 
The results below are stated and proved for
densities that are bounded away from $0$ and $\infty$, uniformly in time, space, and scale. 
In this regime 
many technical arguments are simplified, and this also considerably shortens the presentation.
At the same time, we believe the introduced structure (the cell problem, its dual splitting, the role of the Fisher information, etc.) is sufficiently robust to carry over to the regime without such bounds. 
The non-trivial task of removing the simplifying assumption is left for future work.
In Section~\ref{sec:discussion} we comment on the changes and challenges involved.

\subsection{The Dirichlet-to-Neumann semigroup as a generalized gradient flow}
\label{sec:intro:DtN-GF}

Consider a smooth open bounded connected domain $\Omega$ with boundary
$\Gamma$. We fix positive measures $\pi$ on $\Omega$ and $\omega$ on $\Gamma$,
both absolutely continuous with respect to the corresponding Lebesgue and
Hausdorff measures, and we write $\pi$ and $\omega$ both for the measures and
their densities. We assume that $\log \pi \in C^1(\overline\Omega)$
and $\log\omega\in C(\Gamma)$, which ensures that both densities are
uniformly bounded from above and below on $\Omega$ and $\Gamma$, respectively.

We derive and study the evolution for a pair of measures $(\rho,\gamma)$, where
$\rho$ is a measure on $\Omega$ and $\gamma$ a measure on $\Gamma$, given by
\begin{subequations}\label{eq:system:DtN}
	\begin{alignat}{2}
		0 &=  \div \pra[\Big]{\pi\nabla\frac{\dx\rho}{\dx\pi}} && \text{ in } \Omega  \,; \label{eq:system:DtN:harmonic} \\
		\label{eq:system:DtN:boundary}
		\tr \frac{\dx\rho}{\dx\pi} &= \frac{\dx\gamma}{\dx\omega}  \qquad && \text{ on } \Gamma \,;\\
		\partial_t \gamma &= - \tr \pi \, \partial_n\frac{\dx\rho}{\dx\pi}  \qquad && \text{ on } \Gamma \,. \label{eq:system:DtN:normal}
	\end{alignat}
\end{subequations}
Writing $u \coloneqq \dx\rho/\dx\pi$ for the interior density and $\ub \coloneqq
\dx\gamma/\dx\omega$ for the boundary density, the
equations~\eqref{eq:system:DtN:harmonic} and~\eqref{eq:system:DtN:boundary} are
the $\pi$-weighted Dirichlet problem with boundary datum~$\ub$,
\begin{equation}\label{eq:harmonic:Dirichlet}
	\div(\pi\nabla u) = 0 \text{ in } \Omega\,, \qquad \tr u = \ub \text{ on } \Gamma\,.
 \end{equation}
For $\ub$ in the trace space $H^{1/2}(\Gamma)=\tr H^1(\Omega)$, whose seminorm is~\eqref{eq:def:H12} in Section~\ref{sec:CELL:GS}, this problem has a unique solution $u\in
H^1(\Omega)$, the \emph{$\pi$-harmonic extension} of $\ub$, which we denote
by $u=\calU(\ub)$. Substituting it back into~\eqref{eq:system:DtN:normal} closes
the system as the boundary evolution
\begin{equation}\label{eq:DtN:evolution:gamma}
	\partial_t \gamma = - \tr \pi \, \partial_n \calU\bra*{\frac{\dx\gamma}{\dx\omega}} \qquad \text{on } \Gamma\,,
\end{equation}
in which the right-hand side is the action of the \emph{Dirichlet-to-Neumann operator}
\begin{equation}\label{eq:DtN:operator:strong}
	\DtN : H^{\frac{1}{2}}(\Gamma)\to \dot H^{-\frac{1}{2}}(\Gamma)\,,\qquad \DtN(\ub) \coloneqq  \tr\pi\,   \partial_n \calU(\ub) \,.
\end{equation}
The evolution equation~\eqref{eq:DtN:evolution:gamma} therefore reads
\begin{equation}\label{eq:DtN}
	\partial_t\gamma = -\DtN\bra*{\frac{\dx\gamma}{\dx\omega}}.
\end{equation}
It is shown in~\cite{Hsu86,ArendtElst2011} that $\DtN$ generates a
Markovian semigroup $(e^{-t\DtN})_{t\geq 0}$ on the space of non-negative
densities on~$\Gamma$. We provide a variational characterization of this
evolution and show that it is a generalized gradient flow with respect to the
relative entropy on the boundary $\calH(\gamma | \omega)$.

A first hint of the structure of the gradient flow is given by the entropy
dissipation. For $\ub$ bounded away from $0$ and $\infty$ the divergence theorem
gives
\begin{align*}
    \pderiv{}{t}\calH(\gamma|\omega) &= - \int_\Gamma \log \ub \,  \tr \pi\, \partial_n \calU\bra*{\ub} \dx{\calH^{d-1}} \notag \\
    &= -\int_\Omega \nabla \log \calU(\ub) \cdot \nabla \calU(\ub) \dx\pi = -2 \calI(\calU(\ub))\leq 0 \,,
\end{align*}
where we recognized the Fisher information, defined by
\begin{equation}\label{eq:def:fisher}
    \calI(\rho|\pi) = \frac{1}{2} \int_\Omega \abs{\nabla \log u}^2 \dx\rho
    = \frac{1}{2} \int_\Omega \frac{\abs{\nabla u}^2}{u}\dx\pi
    = 2 \int_\Omega \abs{\nabla\sqrt{u}}^2 \dx\pi \,,
\end{equation}
the first two expressions being valid for $u>0$ and the third for all $u\geq 0$, see
e.g.~\cite[Lemma 2.2]{GianazzaSavareToscani2008}. Since $\calI(\rho|\pi)$
depends on $\rho$ only through its density $u=\dx\rho/\dx\pi$, we abbreviate
\begin{equation}\label{eq:def:fisher:density}
    \calI(u) \coloneqq \calI(u\pi \,|\, \pi)
\end{equation}
and use the two notations interchangeably. The Fisher information satisfies a
minimality property for the $\calU(\ub)\pi$-weighted Dirichlet form and can for
this reason be expressed in terms of a weighted $H^{1/2}$-seminorm of the boundary
datum $\ub$ alone (see Lemma~\ref{lem:logu:optimal}),
\begin{equation}\label{eq:Fisher:H12}
    \calI(\calU(\ub)) = \frac{1}{2} \int_\Omega \abs*{ \nabla \log \calU(\ub)}^2 \calU(\ub) \dx{\pi} = \frac{1}{2} \norm{\log \ub}_{\dot H^{1/2}(\Gamma,\calU(\ub)\pi)}^2 \,.
\end{equation}
We emphasize that, although we use the trace-space notation
$\dot H^{1/2}(\Gamma,\cdot)$, this seminorm depends on the weight
$\calU(\ub)\pi$, which is defined on all of~$\Omega$; see~\eqref{eq:def:H12}.
Identity~\eqref{eq:Fisher:H12} suggests a variational characterization of~\eqref{eq:DtN} obtained by
differentiating the entropy along a curve $\gamma:[0,T]\to\calP_{\mathrm{ac}}(\Gamma)$.
By mass preservation we can parametrize its time derivative by the normal
trace\footnote{For the exact definition of the normal trace we refer to
Definition~\ref{def:limit-CE} and Section~\ref{sec:normaltrace}.} of a divergence-free
flux $J:[0,T]\to\calM(\Omega;\R^d)$, $\div J\equiv0$, such that $\dot\gamma = n\cdot J$
on~$\Gamma$. The Gau\ss--Green theorem and $\div J=0$ then give, for any
sufficiently smooth extension $\calS=\calS(\ub,\dot\gamma)$ of $\ub$ into the interior,
\begin{equation}\label{eq:intro:chain}
    \pderiv{}{t} \calH(\gamma|\omega) = \skp{ \log \ub , \dot \gamma }_\Gamma
        = \int_\Gamma\log \ub \ n \cdot \dx J 
        =  \int_\Omega  \nabla \log \calS \cdot \dx J \,,
\end{equation}
and a weighted Cauchy--Schwarz inequality bounds the right-hand side from below by
\begin{equation}\label{eq:intro:CS}
    - \frac{1}{2} \int_\Omega \frac{1}{\calS}\abs*{\frac{\dx J}{\dx \pi}}^2 \dx\pi
    - \frac{1}{2}\int_\Omega \abs{\nabla \log \calS}^2 \calS \dx\pi \,.
\end{equation}
Equality holds formally if and only if $\dx J = -\nabla\calS\dx\pi$,
which by $\div J=0$ forces~$\calS$ to be $\calU(\ub)$, so that $\dot\gamma = n\cdot J = -\tr\pi\,\partial_n\calU(\ub)$
is again~\eqref{eq:DtN}. In other words, the $\DtN$-evolution is exactly the
equality case in~\eqref{eq:intro:chain}--\eqref{eq:intro:CS}, for \emph{any}
choice of extension~$\calS$ that reduces to $\calU(\ub)$ along solutions.

We choose to resolve this underdetermination by optimizing over the interior
configuration. 
Indeed, fixing $\ub$ and $\dot\gamma=\jb$ on $\Gamma$ and minimizing the
two terms in~\eqref{eq:intro:CS} over all compatible interior states leads to
the \emph{cell formula} (see Definition~\ref{def:cell})
\begin{equation}\label{eq:def:CELL:intro}
    \sfM(\gamma,\jb) \coloneqq \inf_{u,J} 
    \set*{ \frac{1}{2} \int_\Omega \pra*{ \frac{1}{u}\abs*{\frac{\dx J}{\dx \pi}}^2  + \frac{\abs{\nabla u}^2}{u} } \dx{\pi} \;\Bigg|\;
    \begin{array}{@{}c@{}}
        \div J = 0 \text{ in $\Omega$}; \\
        n\cdot J = \jb \text{ on $\Gamma$}; \\
        \tr u = \pderiv{\gamma}{\omega} \text{ on $\Gamma$} ,
    \end{array}}
\end{equation}
and to the variational description of the $\DtN$-evolution~\eqref{eq:DtN} by the
energy-dissipation functional
\begin{equation}\label{eq:DeGiorgi:red}
    \calL^{\DtN}_T(\gamma) = \calH(\gamma_T|\omega) - \calH(\gamma_0|\omega) + \int_0^T \sfM(\gamma_t ,\dot\gamma_t) \dx{t} \,.
\end{equation}
Formally, the passage from~\eqref{eq:intro:CS} to~\eqref{eq:def:CELL:intro} is
an instance of the contraction principle of large-deviation
theory~\cite[Theorem~4.2.1]{DemboZeitouni2010}, i.e., minimizing the interior cost
over all configurations compatible with the prescribed boundary data selects,
among the unlikely interior realizations, the most likely one. This is
consistent with the general correspondence between large-deviation rate
functionals and generalized gradient structures~\cite{MielkePeletierRenger2014}.
However, other ways to resolve the indeterminacy are possible, and we discuss one such way in Section~\ref{sec:discussion}.
\subsection{Main results}
\label{sec:intro:results}

Throughout the paper we fix two constants $0<c\leq C<\infty$ and work with the
convex sets of \emph{nondegenerate densities}
\begin{equation}\label{eq:def:Acc}
    \begin{split}
    \Acc &\coloneqq \set[\big]{u\in H^1(\Omega) \,:\, 
                c\leq u\leq C \ \text{ $\pi$-a.e.\ in }\Omega} , \\
    \Acc_\Gamma &\coloneqq \set[\big]{\ub\in H^{1/2}(\Gamma) \,:\, 
                c\leq \ub\leq C \ \text{ $\omega$-a.e.\ on }\Gamma} .
    \end{split}
\end{equation}
Since $\log\pi\in C^1(\overline\Omega)$, $\log\omega\in C(\Gamma)$, and $c\leq u,\ub\leq C$, the
weighted spaces $H^1(\Omega,u\pi)$, $L^2(\Gamma,\ub\omega)$ for $u\in \Acc, \ub\in\Acc_\Gamma$ coincide with the
unweighted ones with equivalent norms.
Note that $\tr\Acc\subseteq\Acc_\Gamma$, and by the comparison principle one has that $\ub\in\Acc_\Gamma$ implies $\calU(\ub)\in\Acc$.

Our first main result is that the functional~\eqref{eq:DeGiorgi:red} that we postulated above 
arises as the EDP lower limit of a sequence of generalized gradient systems.
These generalized gradient systems describe a diffusion in $\Omega$ that interacts with the boundary as a reservoir, and the limit corresponds to vanishing diffusion in the bulk; see system~\eqref{eq:system:eps}
in Section~\ref{sec:system:eps}.
\begin{svgraybox}\vspace{-0.75\baselineskip}
\begin{result}[EDP lower limit]\label{result:EDPlimit}
    There exists a family of gradient systems
    $\set{(\CE_T,\calE^\eps,\calR^\eps)}_{\eps>0}$ having $\calL^{\DtN}_T$ as
    \emph{EDP lower limit} along nondegenerate families.
    That is, every family of curves satisfying the continuity equation
    with interior and boundary densities in $[c,C]$ and with uniformly bounded dissipations converges, along a
    subsequence, to a limit curve $(\gamma,\jb)$ satisfying the limit
    continuity equation, and along such families with well-prepared initial data we have 
    \begin{equation*}
        \liminf_{\eps\to0}\ \calL^\eps_T(\mu^\eps,\gamma^\eps;J^\eps,j^\eps) \ \geq\ \calL^{\DtN}_T(\gamma) \,.
    \end{equation*}
\end{result}
\vspace{-0.75\baselineskip}\end{svgraybox}
The approximating $\eps$-dependent gradient systems are described in
Section~\ref{sec:system:eps}; the full statement and proof of
Result~\ref{result:EDPlimit} are contained in
Section~\ref{sec:variational:limit}, where the compactness is
Theorem~\ref{thm:comp} and the lower limit inequality is
Theorem~\ref{thm:liminf}. 
Here \emph{well-prepared} means that the initial
energies converge to the limit energy (see~\eqref{hyp:liminf:wellprep}). 

\medskip
Result~\ref{result:EDPlimit} justifies that the reduced EDP
functional~\eqref{eq:DeGiorgi:red} can be realized as a limit. 
It remains to derive from it a gradient system for the $\DtN$-evolution. For this we prove that the cell
formula~\eqref{eq:def:CELL:intro} splits into a primal-dual pair of dissipation
potentials (cf.\ the beginning of Section~\ref{sec:CELL:GS} for the precise definitions of the norms on the right-hand sides below).
\begin{svgraybox}\vspace{-0.75\baselineskip}
\begin{result}[Duality splitting]\label{result:CELL:splitting}
    Let $\gamma = \ub\omega$ with $\ub\in\Acc_\Gamma$.
    There exists a primal-dual pair $(\sfR,\sfR^*)$ of dissipation potentials,
    characterized by the variational problems 
    \begin{align}
        \sfR(\gamma,\jb) &= \inf_{u\in\Acc: \tr u=\ub} 
            \set*{ \frac{1}{2} \norm{\jb}_{\dot H^{-1/2}(\Gamma,u \pi)}^2 
                + \calI(u) - 2\norm{\sqrt\ub}^2_{\dot H^{1/2}(\Gamma,\pi)} } \,, 
                \label{eq:def:R:intro} \\
        \sfR^*(\gamma,\eta) &= \sup_{u\in\Acc: \tr u =\ub} 
            \set*{ \frac{1}{2} \norm{\eta}_{\dot H^{1/2}(\Gamma,u \pi)}^2 
                - \calI(u) + 2\norm{\sqrt\ub}^2_{\dot H^{1/2}(\Gamma,\pi)} } \,,
                \label{eq:def:Rs:intro}
    \end{align}
    such that the cell formula~\eqref{eq:def:CELL:intro} has the \emph{duality splitting}
    \begin{equation}\label{eq:CELL:splitting}
        \sfM(\gamma,\jb) = \sfR(\gamma,\jb) + \sfR^*(\gamma,-\rmD \sfE(\gamma)) \,,
    \end{equation}
    where $\sfE(\gamma) = \calH(\gamma|\omega)$ and $\rmD\sfE(\gamma)=\log\ub$.
    (Here $\rmD$ denotes the derivative with respect to $\gamma$, represented through the measure-function duality pairing $\skp{\gamma,v}_\Gamma\coloneqq \int_\Gamma v\dx\gamma$.)
\end{result}
\vspace{-0.75\baselineskip}\end{svgraybox}
The full statement (Theorem~\ref{thm:splitting}) and its proof are contained in
Section~\ref{sec:CELL:GS}. The proof proceeds by checking that
$\sfR(\gamma,\jb) \coloneqq \sfM(\gamma,\jb)-\sfM(\gamma,0)$ is an admissible
dissipation potential (Lemma~\ref{lem:props:cell}) and
identifying its variational form~\eqref{eq:def:R:intro} and that of its
dual~\eqref{eq:def:Rs:intro} (Lemma~\ref{lem:M:R:Rstar:char}). In particular the
$\DtN$-evolution~\eqref{eq:DtN} carries the generalized gradient structure
induced by $(\sfE,\sfR)$,
\begin{equation}\label{eq:DtN:GS}
    \dot \gamma = \rmD_2 \sfR^*\bigl(\gamma,-\rmD \sfE(\gamma)\bigr) \,,
\end{equation}
with $\rmD_2$ the differential in the second argument.

The identity~\eqref{eq:CELL:splitting} also shows that the dual 
dissipation potential evaluated at the entropic driving force is the \emph{no-flux} value of the cell formula,
\[
    \sfR^*\bigl(\gamma,-\rmD\sfE(\gamma)\bigr) = \sfM(\gamma,0) = 2\norm{\sqrt\ub}^2_{\dot H^{1/2}(\Gamma,\pi)} = \calI\bigl(\calU(\sqrt\ub)^2\bigr) \quad\text{with }\gamma = \ub \omega\,,
\]
that is, the Fisher information of the \emph{Fisher-harmonic} extension, the
minimizer of $u\mapsto\calI(u)$ at fixed boundary datum. The primal potential
$\sfR$ carries the remaining, flux-dependent part. The two extensions occurring
here are genuinely different: $\calU(\ub)$ minimizes the Dirichlet energy, whereas $\calU(\sqrt\ub)^2$ minimizes the Fisher information and appears in the value $\sfM(\gamma,0)$ (see Lemma~\ref{lem:props:cell}).

\begin{remark}[Full gradient structure]
For the claim that the EDP functional~\eqref{eq:DeGiorgi:red} defines a gradient structure for the $\DtN$ evolution one would normally require it to be non-negative and for the optimal flux to satisfy a suitable constitutive relation.
This amounts to proving a chain rule (see e.g.~\cite{RossiSavare06,MielkeRossiSavare13,Mielke16a}). 
We leave this question to future work.
\end{remark}

\section{The approximating diffusion problems}\label{sec:system:eps}

Let $\eps >0$ and consider the coupled bulk--boundary problem
\begin{subequations}\label{eq:system:eps}
\begin{alignat}{2}
\partial_t \mu 
&= \eps^{-1}\div \pra[\Big]{\mu\nabla\log\bra[\Big]{\frac1\eps\frac{\dx\mu}{\dx\pi}}} \qquad&& \text{ in } \Omega\,; \\
- \eps^{-1} \tr \mu \,
\partial_n \log\bra[\Big]{\frac1\eps\frac{\dx\mu}{\dx\pi}}
&= \eps^{-\tfrac{1}{2}}  \bra*{ \tr\!\left(\frac1\eps\frac{\dx\mu}{\dx\pi}\right) - \frac{\dx\gamma}{\dx\omega}}\sqrt{\tr\pi\,\omega} \qquad && \text{ on } \Gamma ; \\
\partial_t \gamma 
&= \eps^{-\tfrac{1}{2}}  \bra*{ \tr\!\left(\frac1\eps\frac{\dx\mu}{\dx\pi}\right) - \frac{\dx\gamma}{\dx\omega}}\sqrt{\tr\pi\,\omega}   \qquad && \text{ on } \Gamma \,.
\end{alignat}
\end{subequations}
The interior measure carries a factor $\eps$, and the natural object
for taking 
the limit is the rescaled measure $\rho^\eps \coloneqq \mu/\eps$ with density $u^\eps \coloneqq
\dx\rho^\eps/\dx\pi$, in terms of which~\eqref{eq:system:eps} becomes
\begin{subequations}\label{eq:system:eps_rescaled}
\begin{alignat}{2}
\label{eq:system:eps_rescaled:interior}
\eps\partial_t \rho^\eps 
&= \div \pra[\big]{\rho^\eps\nabla\log u^\eps} && \text{ in } \Omega\,; \\
\label{eq:system:eps_rescaled:Robin}
- \tr \rho^\eps \,\partial_n \log u^\eps
&= \eps^{-\tfrac{1}{2}}  \bra*{ \tr u^\eps - \ub}\sqrt{\tr\pi\,\omega} \qquad && \text{ on } \Gamma ; \\
\label{eq:system:eps_rescaled:boundary}
\partial_t \gamma 
&= \eps^{-\tfrac{1}{2}}  \bra*{ \tr u^\eps - \ub}\sqrt{\tr\pi\,\omega}   \qquad && \text{ on } \Gamma \,.
\end{alignat}
\end{subequations}
The two sides in~\eqref{eq:system:eps_rescaled:Robin} scale differently, and the
left-hand side of~\eqref{eq:system:eps_rescaled:interior} formally vanishes as
$\eps\to0$, while the formal limit of the right-hand side of~\eqref{eq:system:eps_rescaled:Robin} enforces the trace identity.
Hence, plugging~\eqref{eq:system:eps_rescaled:Robin} into~\eqref{eq:system:eps_rescaled:boundary}, we expect the limit system to be~\eqref{eq:system:DtN} of
Section~\ref{sec:intro:DtN-GF}.

Before turning to the gradient structure, we record that for strong solutions of~\eqref{eq:system:eps_rescaled} bounds of the form~\eqref{eq:def:Acc} are propagated in time. 
To this end, note that written in the
densities $(u^\eps,\ub)$, and abbreviating $\beta \coloneqq \sqrt{\tr\pi\,\omega}>0$
for the exchange weight on $\Gamma$, the
system~\eqref{eq:system:eps_rescaled} reads
\begin{subequations}\label{eq:system:eps:density}
\begin{alignat}{2}
\label{eq:system:eps:density:interior}
\eps\,\pi\,\partial_t u^\eps &= \div\bra[\big]{\pi\nabla u^\eps} && \text{ in } \Omega\,; \\
\label{eq:system:eps:density:Robin}
- \tr\pi\,\partial_n u^\eps &= \eps^{-\tfrac12}\,\beta\,\bra[\big]{\tr u^\eps - \ub} \qquad && \text{ on } \Gamma\,; \\
\label{eq:system:eps:density:ODE}
\omega\,\partial_t \ub &= \eps^{-\tfrac12}\,\beta\,\bra[\big]{\tr u^\eps - \ub} && \text{ on } \Gamma\,.
\end{alignat}
\end{subequations}
The exchange term enters the interior flux~\eqref{eq:system:eps:density:Robin}
and the boundary evolution~\eqref{eq:system:eps:density:ODE} with opposite
signs, and it is monotone in $\tr u^\eps-\ub$. 
This allows us to prove the following lemma.
\begin{lemma}[Propagation of two-sided bounds]\label{lem:comparison}
Let $\eps>0$ and let $(u^\eps,\ub)$ be a classical solution
of~\eqref{eq:system:eps:density} on $[0,T]$, smooth up to the boundary. If
\[
    c\leq u^\eps_0\leq C \ \text{ in }\Omega
    \qquad\text{and}\qquad
    c\leq \ub_0\leq C \ \text{ on }\Gamma
\]
for constants $0<c\leq C<\infty$, then the same bounds hold at every
$t\in[0,T]$,
\[
    c\leq u^\eps_t\leq C \ \text{ in }\Omega
    \qquad\text{and}\qquad
    c\leq \ub_t\leq C \ \text{ on }\Gamma \,.
\]
\end{lemma}
\begin{proof}
Throughout the proof we fix and drop $\eps>0$, which does not change any of the arguments.
The system~\eqref{eq:system:eps:density} is linear and every pair of constants $(k,k)$ with $k\in\R$ solves it.
Hence $(u-C,\ub-C)$ and $(c-u,c-\ub)$ are again solutions, and both are non-positive at time zero.
It is therefore sufficient to show that solutions of~\eqref{eq:system:eps:density} preserve non-positivity.
Applied to the first pair this gives the upper bound, applied to the second pair the lower bound.

We therefore assume that $u$ and $\ub$ are solutions of~\eqref{eq:system:eps:density} that are non-positive at time zero, and set
\[
    \phi_t \coloneqq u_t^+ \ \text{ on }\Omega \,,
    \qquad
    \psi_t \coloneqq \ub_t^+ \ \text{ on }\Gamma \,.
\]
These vanish at $t=0$ by assumption, and  satisfy $\tr\phi_t = (\tr u_t)^+$ since $u_t$ is continuous up to $\Gamma$. Multiplying the interior
equation~\eqref{eq:system:eps:density:interior} by $\phi_t$, integrating over
$\Omega$ and using the Gau\ss--Green formula together with the pointwise
identities $\nabla u_t\cdot\nabla\phi_t = \abs{\nabla\phi_t}^2$ and
$\partial_t u_t\,\phi_t = \tfrac12\partial_t\phi_t^2$ gives
\[
    \frac{\eps}{2}\pderiv{}{t}\int_\Omega\phi_t^2\dx\pi
    = -\int_\Omega\abs{\nabla\phi_t}^2\dx\pi
      + \int_\Gamma \tr\pi\,\partial_n u_t\;(\tr u_t)^+\dx\calH^{d-1} \,.
\]
Multiplying the boundary equation~\eqref{eq:system:eps:density:ODE} by $\psi_t$
and integrating over $\Gamma$ gives in the same way
\[
    \frac12\pderiv{}{t}\int_\Gamma\psi_t^2\dx\omega
    = \frac{1}{\sqrt\eps}\int_\Gamma \beta\,(\tr u_t - \ub_t)\,(\ub_t)^+\dx\calH^{d-1} \,.
\]
Inserting the boundary condition~\eqref{eq:system:eps:density:Robin} into the
first identity and adding the two, we find
\begin{align}\label{eq:comparison:energy}
    \pderiv{}{t} &\Bigl[\frac{\eps}{2}\int_\Omega\phi_t^2\dx\pi + \frac12\int_\Gamma\psi_t^2\dx\omega\Bigr]\\
    &= -\int_\Omega\abs{\nabla\phi_t}^2\dx\pi
      - \frac{1}{\sqrt\eps}\int_\Gamma \beta\,(\tr u_t-\ub_t)\bigl[(\tr u_t)^+ - \ub_t^+\bigr]\dx\calH^{d-1} \,. \notag
\end{align}
Both terms on the right-hand side are non-positive. 
Indeed, for the first this is immediate, and
for the second, it follows since $\beta\geq0$ and $s\mapsto s^+$ is non-decreasing. 
Hence the bracket on the left of~\eqref{eq:comparison:energy} is non-increasing,
vanishes at $t=0$, and is non-negative.
Thus, it vanishes on $[0,T]$, i.e., $u_t\leq 0$ in $\Omega$ and $\ub_t\leq 0$ on $\Gamma$.
\end{proof}

We now introduce the prelimit gradient system for $\eps>0$ fixed.
Denote by $\pi^\eps \coloneqq \eps\,\pi$ the rescaled interior reference measure. The
system~\eqref{eq:system:eps} is a gradient flow driven by the combined relative
entropy
\begin{equation}\label{eq:def:combined:entropy}
\calE^\eps(\mu,\gamma) = \calH(\mu|\pi^\eps) + \calH(\gamma|\omega) = \eps\,\calH(\rho^\eps|\pi) + \calH(\gamma|\omega)\,,
\end{equation}
with respect to a generalized gradient system in continuity equation
format (cf.\ \cite{PeletierSchlichting2023}), which we now specify. 
The formal continuity equation is
\begin{equation}\label{eq:formal:CE}
    \partial_t \mu_t + \div J_t = 0 \text{ in } \Omega,
    \qquad
    \partial_t \gamma_t  = j_t = J_t\cdot n  \text{ on }\Gamma,
\end{equation}
to be understood in the sense of distributions. 
Similar to~\cite[Def.~3.5]{HeinzePietschmannSchlichting2024},
the normal flux in~\eqref{eq:formal:CE} produces a discrete bulk-interface gradient.
We now specify this and rigorously define the continuity equation.
\begin{definition}[Joint continuity equation and bulk-interface gradient]\label{def:CE}
The \emph{bulk-interface gradient} of a test pair
$(\varphi,\phi)\in C^1(\overline\Omega)\times C^1(\Gamma)$ is
\begin{equation}\label{eq:def:anabla}
    \anabla(\varphi,\phi) \coloneqq \bigl(\nabla\varphi,\,\phi - \tr\varphi\bigr) \in C(\overline\Omega;\R^d) \times C(\Gamma)\,,
\end{equation}
with the dual pairing against a flux pair $(J,j) \in \calM(\Omega;\R^d)\times \calM(\Gamma)$,
\begin{equation}\label{eq:def:anabla:pairing}
    \skp[\big]{(J,j),\anabla(\varphi,\phi)}_{\Omega,\Gamma}
    \coloneqq \int_\Omega \nabla\varphi \cdot \dx J + \int_\Gamma \bigl(\phi - \tr\varphi\bigr)\dx j\,.
\end{equation}
A quadruple
$(\mu,\gamma,J,j)$ with $(\mu,\gamma)\in \mathrm{AC}([0,T];\calM_{\geq0}(\Omega))\times\mathrm{AC}([0,T];\calM_{\geq0}(\Gamma))$ and $(J,j)\in(\calM([0,T]\times\Omega;\R^d) \times\calM([0,T]\times\Gamma))$
satisfies the \emph{joint continuity equation},
$(\mu,\gamma,J,j)\in\CE_T$, if for all
$(\varphi,\phi)\in C^1([0,T]\times\overline\Omega)\times C^1([0,T]\times\Gamma)$
\begin{align}\label{eq:CE:weak}
        \MoveEqLeft \int_{\Omega}\varphi(T,x)\dx\mu_T + \int_\Gamma\phi(T,x)\dx\gamma_T
        = \int_{\Omega} \varphi(0,x) \dx{\mu_0}  + \int_{\Gamma} \phi(0,x) \dx{\gamma_0} \\
        &+ \int_0^T \bra[\bigg]{ \int_\Omega \partial_t\varphi \dx{\mu_t} + \int_\Gamma \partial_t\phi \dx{\gamma_t} + \skp[\big]{(J_t,j_t),\anabla(\varphi,\phi)}_{\Omega,\Gamma} } \dx{t} \,. \notag
    \end{align}
\end{definition}
Note that the weak form~\eqref{eq:CE:weak} encodes the normal flux condition $J\cdot n = j$.
Indeed, for smooth fluxes (or fluxes with measure-valued divergence, cf.~\eqref{eq:Gauss-Green} below) an integration by parts in~\eqref{eq:def:anabla:pairing} gives rise to the additional boundary term
$\int_\Gamma \tr\varphi\,\bigl(n\cdot \dx J - \dx j\bigr)$, which vanishes for all $\varphi$ if and only if $J\cdot n=j$.
Furthermore, the narrow absolute continuity in time of the curve $(\mu,\gamma)$ is already implied by the integrability of the fluxes $\int_0^T\bigl(\abs{J_t}(\Omega) + \abs{j_t}(\Gamma)\bigr)\dx t <\infty$. 
\begin{definition}[(Dual) dissipation potentials and relaxed slopes]\label{def:potentials}
For a bulk force $\Xi\in C(\overline\Omega;\R^d)$ and a boundary forcing
$\zeta\in C(\Gamma)$ we set
\begin{subequations}\label{eq:dissipation:dual}
    \begin{align}
    (\calR^\eps)^*\bigl((\mu,\gamma),(\Xi,\zeta)\bigr) &\coloneqq (\calR_\Omega^\eps)^*(\mu,\Xi) + (\calR_{\Gamma}^\eps)^*\bigl((\tr\mu,\gamma),\zeta\bigr)\,,
    \label{eqdef:Rstar-both}\\
    (\calR_\Omega^\eps)^*(\mu,\Xi) &\coloneqq \frac{1}{2\eps} \int_{\Omega} \abs*{\Xi}^2 \dx\mu = \frac{1}{2} \int_{\Omega} \abs*{\Xi}^2 \dx\rho^\eps \,, 
    \label{eqdef:Rstar-Omega}\\
    (\calR_{\Gamma}^\eps)^*\bigl((\tr\mu,\gamma),\zeta\bigr) 
    &\coloneqq \frac{1}{\eps} \int_{\Gamma} \sfC^*(\zeta) \dx\sqrt{\tr\mu\,\gamma}
    = \frac{1}{\sqrt{\eps}} \int_{\Gamma} \sfC^*(\zeta) \dx\sqrt{\tr\rho^\eps\,\gamma} \,,
    \label{eqdef:Rstar-Gamma}
\end{align}
\end{subequations}
where $\sfC^*(\zeta)\coloneqq 4\bra*{\cosh(\zeta/2)-1}$, whose convex conjugate is
\begin{equation*}%
    \sfC(s) \coloneqq \sup_{\zeta\in\R}\set*{s\zeta - \sfC^*(\zeta)} = 2s\,\arsinh\bra*{\frac{s}{2}} - 2\sqrt{4+s^2} + 4 \,.
\end{equation*}
The (primal) dissipation potentials are the convex conjugates of the dual dissipation potentials w.r.t.\ the second argument. 
They are given by
\begin{subequations}\label{eq:dissipation:primal}
\begin{align}
\calR^\eps\bigl((\mu,\gamma),(J,j)\bigr) &= \calR_\Omega^\eps(\mu,J) + \calR_\Gamma^\eps\bigl((\tr\mu,\gamma),j\bigr)\,,\\
\calR_\Omega^\eps(\mu,J) &= \frac{1}{2} \int_\Omega \abs*{\frac{\dx J}{\dx\rho^\eps}}^2 \dx\rho^\eps\,, \label{eqdef:R:Omega}\\
\calR_\Gamma^\eps\bigl((\tr\mu,\gamma),j\bigr) &=\frac1{\sqrt\eps} \int_\Gamma \sfC\bra*{\frac{\sqrt{\eps}\,\dx j}{\dx\sqrt{\tr\rho^\eps\,\gamma}}}\dx\sqrt{\tr\rho^\eps\,\gamma}\,.
\end{align}
\end{subequations}
Finally, we define the relaxed slopes
in terms of $u^\eps = \pderiv{\mu}{{\pi^\eps}} = \pderiv{\rho^\eps}{\pi}$ and $\ub = \pderiv{\gamma}{\omega}$ by
\begin{subequations}\label{eq:slopes}
\begin{align}
    \calS^\eps(\mu,\gamma) &= \calS_{\Omega}^\eps(\mu) + \calS_{\Gamma}^\eps\bigl((\tr\mu,\gamma)\bigr)\,,\\
    \calS_{\Omega}^\eps(\mu) &= 
    2 \int_\Omega \abs*{\nabla \sqrt{u^\eps}}^2 \dx\pi = \calI(u^\eps)\,, \label{eqdef:Slope:Omega} \\
    \calS_{\Gamma}^\eps\bigl((\tr\mu,\gamma)\bigr) 
    &= \frac{2}{\sqrt{\eps}} \int_{\Gamma} \abs[\big]{\sqrt{\tr u^\eps} - \sqrt{\ub}}^2 \dx\sqrt{\tr\pi\,\omega}\,.
\end{align}
\end{subequations}
with $\calI$ the Fisher information~\eqref{eq:def:fisher}.
\end{definition}
The interior potential~\eqref{eqdef:R:Omega} and the interior slope~\eqref{eqdef:Slope:Omega} are independent of~$\eps$ in the
rescaled variable $\rho^\eps$ (resp.\ its density $u^\eps$) and we define
\begin{equation}\label{eqdef:RR*:cell}
    \calR_\Omega^\ast(\rho,\Xi) \coloneqq \frac{1}{2} \int_{\Omega} \abs*{\Xi}^2 \dx\rho\,
    \quad\text{and}\quad
    \calR_\Omega(\rho,J)  \coloneqq  \frac{1}{2} \int_\Omega \abs*{\frac{\dx J}{\dx\rho}}^2 \dx\rho \,,
\end{equation}
with the convention that $\calR_\Omega(\rho,J)=+\infty$ if $J\not\ll \rho$. 
Together with the Fisher information $\calI$, they are the central objects for the cell formula, below.

In contrast, the boundary contributions carry an
$\eps^{\pm1/2}$-dependence reflecting the fast exchange of mass between $\Omega$
and~$\Gamma$, which produces the boundary
condition~\eqref{eq:system:DtN:boundary} in the limit.

\begin{definition}[Energy-dissipation functional]\label{def:Leps}
For $(\mu,\gamma;J,j)\in \CE_T$ we set
\begin{align}
\calL^\eps_T(\mu,\gamma;J,j) &= \calE^\eps(\mu_T,\gamma_T) - \calE^\eps(\mu_0,\gamma_0) + \calD^\eps_T(\mu,\gamma;J,j) \,, \\
    \label{eqdef:Deps}
\calD^\eps_T(\mu,\gamma;J,j) &=  \int_0^T \bigl[ \calR_{\Omega}^\eps(\mu,J) + \calR_{\Gamma}^\eps\bigl((\tr\mu,\gamma),j \bigr) + \calS_{\Omega}^\eps(\mu) + \calS_{\Gamma}^\eps\bigl((\tr\mu,\gamma)\bigr) \bigr] \dx{t} \,.
\end{align}
A curve with $\calL^\eps_T=0$ is called an \emph{EDE-solution} of the gradient system
$(\CE_T,\calE^\eps,\calR^\eps)$. 
A curve with $\calL^\eps_T\leq0$ is called an
\emph{EDI-solution}.
\end{definition}
We emphasize that due to the superlinearity of the square and $\sfC$, respectively, for each fixed $\eps>0$ the boundedness of the dissipation, $\calD^\eps_T(\mu,\gamma;J,j)<\infty$, implies the time integrability of the fluxes, $\int_0^T\bigl(\abs{J_t}(\Omega) + \abs{j_t}(\Gamma)\bigr)\dx t <\infty$.

For $(\mu,\gamma)$ smooth and strictly positive the entropic driving force is
the pair $\rmD\calE^\eps(\mu,\gamma) = (\log(\dx\mu/\dx\pi^\eps),\,
\log(\dx\gamma/\dx\omega))$, and applying $\anabla$ component-wise turns it into
an admissible argument for $(\calR^\eps)^*$,
\begin{equation}\label{eq:anabla:DE}
    \anabla \rmD\calE^\eps(\mu,\gamma) = \biggl(\nabla\log\frac{\dx\mu}{\dx\pi^\eps},\ \log\frac{\dx\gamma}{\dx\omega} - \log\tr\frac{\dx\mu}{\dx\pi^\eps}\biggr) \,.
\end{equation}
Furthermore, in this case we have $\calS^\eps(\mu,\gamma) =
(\calR^\eps)^*((\mu,\gamma),-\anabla\rmD\calE^\eps(\mu,\gamma))$, and the chain
rule reads $\pderiv{}{t}\calE^\eps(\mu_t,\gamma_t) = \skp{(J_t,j_t),\anabla
\rmD\calE^\eps(\mu_t,\gamma_t)}_{\Omega,\Gamma}$. 
Combining it with the Fenchel--Young inequality yields $\calL^\eps_T\geq0$ for every curve in $\CE_T$,
with equality if and only if the constitutive relation
\begin{equation}\label{eq:contact:eps}
    (J_t,j_t) = \rmD_2 (\calR^\eps)^*\bigl((\mu_t,\gamma_t),\,-\anabla\rmD\calE^\eps(\mu_t,\gamma_t)\bigr) 
\end{equation}
holds for a.e.\ $t$, which is precisely the strong
form~\eqref{eq:system:eps}. 

In particular, if $(\mu_t,\gamma_t)_{t\in[0,T]}$ is a smooth, strictly
positive classical solution of~\eqref{eq:system:eps}, with fluxes
$J_t \coloneqq -\rho^\eps_t\nabla\log u^\eps_t$ and
$j_t \coloneqq \eps^{-1/2}(\tr u^\eps_t - \ub_t)\sqrt{\tr\pi\,\omega}$, then the derivative of the energy is given by
\begin{equation}\label{eq:EDR:eps}
\begin{split}
\pderiv{}{t}\calE^\eps(\mu_t,\gamma_t) 
= &-\int_\Omega u^\eps_t \abs{\nabla\log u^\eps_t}^2 \dx\pi \\
&- \frac{1}{\sqrt{\eps}}\int_\Gamma \bigl(\tr u^\eps_t - \ub_t\bigr)\bigl(\log \tr u^\eps_t - \log \ub_t\bigr)\dx\sqrt{\tr\pi\,\omega}\,.
\end{split}
\end{equation}
Note that both integrals are non-negative, and after integration in time we obtain the energy-dissipation identity
\begin{equation}\label{eq:EDB:eps:concrete}
    \calE^\eps(\mu_T,\gamma_T) + \int_0^T\bigl[\calR^\eps(\mu_t,\gamma_t;J_t,j_t) + \calS^\eps(\mu_t,\gamma_t)\bigr]\dx t = \calE^\eps(\mu_0,\gamma_0)\,.
\end{equation}
That is, $(\mu,\gamma;J,j)$ is an EDE-solution.

Indeed, 
we calculate the derivative of the interior entropy,
which by the chain rule and integration by parts takes the form
\begin{align*}
    \pderiv{}{t}\calH(\mu_t|\pi^\eps)
    &= \frac{1}{\eps}\int_\Omega \log u^\eps_t\, \div\bigl[\mu_t\nabla\log u^\eps_t\bigr]\dx x \\
    &= -\int_\Omega u^\eps_t\abs{\nabla\log u^\eps_t}^2\dx\pi + \frac{1}{\sqrt\eps}\int_\Gamma \log \tr u^\eps_t\,(\ub_t-\tr u^\eps_t)\dx\sqrt{\tr\pi\,\omega} \,,
\end{align*}
using $\mu=\eps u^\eps\pi$ and the boundary condition
$\tfrac{1}{\eps}\tr[\mu\nabla\log u^\eps]\cdot n = \tfrac{1}{\sqrt\eps}(\ub - \tr
u^\eps)\sqrt{\tr\pi\,\omega}$.
The derivative of the boundary entropy amounts to
$\pderiv{}{t}\calH(\gamma_t|\omega) = \tfrac{1}{\sqrt\eps}\int_\Gamma
\log\ub_t\,(\tr u^\eps_t-\ub_t)\dx\sqrt{\tr\pi\,\omega}$. Adding the two
yields~\eqref{eq:EDR:eps}. The identification with $\calR^\eps+\calS^\eps$
in~\eqref{eq:EDB:eps:concrete} is a direct computation, carried out, e.g., in~\cite{HeinzePietschmannSchlichting2024} for a one-dimensional domain.
For the gradient-flow flux given by~\eqref{eq:contact:eps},
one finds $\calR^\eps_\Omega = \calS^\eps_\Omega =
\tfrac12\int_\Omega u^\eps|\nabla\log u^\eps|^2\dx\pi$, while on the boundary
$\sqrt\eps\,\dx j/\dx\nu = 2\sinh(\zeta/2)$ for $\zeta\coloneqq\log(\tr u^\eps/\ub)$ and
$\nu\coloneqq\sqrt{\tr\rho^\eps\gamma}$, so that the
identity $\sfC\bigl(2\sinh(\zeta/2)\bigr) + \sfC^*(\zeta) = 2\zeta\sinh(\zeta/2)$ turns $\calR^\eps_\Gamma+\calS^\eps_\Gamma$
into the second integral in~\eqref{eq:EDR:eps}.

\section{Compactness and the lower limit cell formula}\label{sec:variational:limit}

In this section we introduce the limit continuity equation and the cell problem,
and we prove that curves of the $\eps$-problem converge to curves of the limit
problem with a lower limit inequality for the energy-dissipation functionals.
This is Result~\ref{result:EDPlimit}.

All measures occurring below are finite, and for measures we use throughout
\emph{narrow} convergence, that is, testing against bounded continuous
functions: $\varrho_n\to\varrho$ narrowly in $\calM(\Omega)$ if and only if
$\int\varphi\dx\varrho_n\to\int\varphi\dx\varrho$ for every $\varphi\in
C_b(\Omega)$. 
On the compact boundary $\Gamma$ narrow convergence is the same as weak-$^\ast$
convergence in $\calM(\Gamma)=C(\Gamma)^\ast$. 

\subsection{The normal trace of divergence-free measure fields}
\label{sec:normaltrace}

The limit continuity equation below transports the boundary measure by a nonlocal flux on $\Gamma$, which is obtained as the normal trace from a divergence-free flux in $\Omega$.
We make those definitions rigorous based on the divergence-measure fields~\cite{ChenFrid1999,ChenFrid2003}.

Let $J\in\calM(\Omega;\R^d)$ be a finite measure field with $\div J=0$ in the
sense of distributions. Then, $J$ is an extended divergence-measure field,
$J\in\calD\calM^{\mathrm{ext}}(\Omega;\R^d)$, in the sense
of~\cite{ChenFrid2003}. By~\cite[Theorem~3.1]{ChenFrid2003}, such a field
admits a \emph{normal trace}\footnote{The actual space used in~\cite[Theorem~3.1]{ChenFrid2003} is $C^{1,\delta}(\Gamma)$ for every $\delta>0$, i.e., differentiable functions with $\delta$-Hölder continuous derivative.}
$\jb\in\bigl(C^{\infty}(\Gamma)\bigr)^\ast$, characterized by the
Gau\ss--Green formula
\begin{align}\label{eq:Gauss-Green}
    \skp{\jb,\tr\varphi}_\Gamma \coloneqq \int_\Omega \varphi\dx(\div J) + \int_\Omega \nabla \varphi \cdot \dx J
    \qquad \text{for all } \varphi\in C^{\infty}(\overline\Omega)\,.
\end{align}
We will sometimes write $\jb = J\cdot n$ and call $\jb$ the boundary flux of $J$.
Observe that the right-hand side of~\eqref{eq:Gauss-Green} vanishes for test functions that are identically zero on $\Gamma$.
Hence, by linearity, $\jb$ is indeed a functional on the boundary data $\tr\varphi$ alone.

Two consequences of~\eqref{eq:Gauss-Green} are used repeatedly.
First, testing with $\varphi\equiv1$ and using $\div J=0$ gives $\skp{\jb,1}_\Gamma=0$, i.e., $\jb$ annihilates constants. 
Second, testing with a general $\varphi$ and using $\div J=0$ gives the bound
\begin{equation}\label{eq:jb:bound}
    \abs{\skp{\jb,\tr\varphi}_\Gamma} \leq \norm{\nabla\varphi}_{L^\infty}\abs{J}(\Omega) \,.
\end{equation}

Throughout the section we apply this to a time-dependent field, that is, to
$J_t$ for almost every $t$, and write $\jb_t$ for the corresponding normal
trace. This is also how the boundary flux arises in Theorem~\ref{thm:comp}. It
is obtained from the limit $J$ of the interior fluxes $(J^\eps)_\eps$, and not
from the traces of the $J^\eps$ themselves, which need not exist in the prelimit.

\subsection{The limit continuity equation and the cell problem}
\label{sec:limitCE}

Before giving the definition, we indicate how the continuity equation emerges in the limit.
The joint continuity equation~\eqref{eq:CE:weak} couples the interior and the boundary
through the exchange term in~\eqref{eq:def:anabla}. 
Choosing a single test function, $\phi=\tr\varphi$ with $\varphi\in C_c^\infty((0,T)\times\overline\Omega)$, this term vanishes and~\eqref{eq:CE:weak} reduces to
\begin{equation}\label{eq:CE:embedded}
    0 = \int_0^T\!\!\int_\Gamma\partial_t\varphi\dx\gamma^\eps + \int_0^T\!\!\int_\Omega\partial_t\varphi\dx\mu^\eps + \int_0^T\!\!\int_\Omega\nabla\varphi\cdot\dx J^\eps \,,
\end{equation}
which is an identity that no longer involves the exchange flux $j^\eps$. 
We will see that in the limit $\eps\to0$ the interior mass disappears, so that only the first and the third
term remain, i.e., the change of the limiting boundary measure $\gamma$ is determined nonlocally by the normal trace of the limit interior flux $J$, which by the discussion in Section~\ref{sec:normaltrace} is well defined.
\begin{definition}[Limit continuity equation]\label{def:limit-CE}
    Let $\gamma\in\mathrm{AC}([0,T];\calM_{\geq0}(\Gamma))$ with $\gamma_t\ll\omega$,
    and let $\jb=(\jb_t)_{t\in[0,T]}$ be a family in $\bigl(C^\infty(\Gamma)\bigr)^\ast$.
    We say that $(\gamma,\jb)$ satisfies the \emph{limit continuity equation}, in short
    $(\gamma,\jb)\in\lCE_T$, if the following two conditions hold.
    \begin{enumerate}[label=(\roman*)]
        \item There exists $J\in\calM([0,T]\times\Omega;\R^d)$ with
        $\abs{J}(\cdot\times\Omega)\ll\scrL^1$ and disintegration
        $J(\dx t,\dx x)=J_t(\dx x)\dx t$ such that, for a.e.\ $t\in[0,T]$,
        $\div J_t=0$ and $\jb_t=J_t\cdot n$ in the sense of~\eqref{eq:Gauss-Green}.
        We call any such $J$ (resp.\ $J_t$) an \emph{interior flux associated with} $\jb$ (resp.\ $\jb_t$).
        \item For all $\phi\in C^\infty([0,T]\times\Gamma)$,
        \begin{equation}\label{eq:limit-CE}
            \int_\Gamma \phi(T,\cdot) \dx\gamma_T - \int_\Gamma \phi(0,\cdot) \dx\gamma_0
            = \int_0^T \!\!\int_\Gamma \partial_t\phi \dx\gamma_t\dx t
            + \int_0^T\skp{\jb_t,\phi(t,\cdot)}_\Gamma\dx t\,.
        \end{equation}
    \end{enumerate}
\end{definition}
Note that since for a.e.\ $t\in[0,T]$ the boundary flux $\jb_t$ is the normal trace of some interior flux $J_t$, condition~\eqref{eq:limit-CE} is equivalent to
\begin{align}\label{eq:limit-CE:J}
    \int_\Gamma \tr\varphi(T,\cdot) \dx\gamma_T - \int_\Gamma \tr\varphi(0,\cdot) \dx\gamma_0 &= \int_0^T \!\!\int_\Gamma \partial_t\tr\varphi\dx\gamma_t\dx t + \int_0^T\!\!\int_\Omega \nabla\varphi\cdot \dx J_t\dx t\,,
\end{align}
for all $\varphi\in C^\infty([0,T]\times\overline\Omega)$.
However, we highlight that this interior flux $J$ is not unique, since adding a divergence-free field with vanishing normal trace does not change $\jb$.

We now introduce the time-independent cell problem, which will later be considered at a.e.\ $t\in[0,T]$
for a pair $(\gamma_t,\jb_t)$.
\begin{definition}[Cell problem]\label{def:cell}
For $(\gamma,\jb)\in \calM_{\geq0}(\Gamma)\times(C^\infty(\Gamma))^\ast$, we introduce the \emph{cell problem}
\begin{equation}\label{eq:CELL}
         \sfM(\gamma,\jb) \coloneqq  \inf_{(u,J)\in M(\gamma,\jb)}
        \Bigl\{ \calR_\Omega\bigl(u\pi, J \bigr) + \calI(u) \Bigr\} \,,
\end{equation}
with $\sfM(\gamma,\jb) \coloneqq +\infty$ if $M(\gamma,\jb)=\emptyset$, where
$M(\gamma,\jb)$ consists of all pairs $(u,J)\in\Acc\times\calM(\Omega;\R^d)$
such that
\begin{enumerate}[label=(\roman*)]
    \item\label{def:cell:i} $J$ is an interior flux associated with $\jb$, that is $\div J \equiv 0$ in the sense of distributions and
    \begin{align}\label{eq:cell:flux}
        \int_\Omega \nabla \varphi \cdot \dx J = \skp{\jb,\tr\varphi}_\Gamma
        \qquad\text{for all }\varphi\in C^\infty(\overline\Omega)\,;
    \end{align}
    \item\label{def:cell:ii} $u$ satisfies the trace condition
    \begin{equation}\label{eq:u-gamma_coupling}
            \tr u = \frac{\dx \gamma}{\dx\omega} \,.
    \end{equation}
\end{enumerate}
\end{definition}

Note that $\sfM(\gamma,\jb)=+\infty$ unless $\gamma=\ub\omega$ with
$\ub\in\tr\Acc\subseteq\Acc_\Gamma$.
Conversely, for $\ub\in\Acc_\Gamma$ the
set $M(\gamma,\jb)$ is non-empty as soon as an admissible $J$ exists, since
$\calU(\ub)\in\Acc$ has $\tr\calU(\ub)=\ub$. 
Within $\Acc$ every competitor satisfies $\log u\in
H^1(\Omega,u\pi)$ and $\tr\log u=\log\tr u$.

For the convenience of the reader, we recall the lower semicontinuity of the two ingredients of the cell formula (also see e.g.~\cite[Lemma 2.2]{GianazzaSavareToscani2008}).
\begin{lemma}[Lower semicontinuity of $\calR_\Omega$ and $\calI$]\label{lem:FI:lower:semicont}
The Fisher information $\calI: \calM_{\geq0}(\Omega) \to [0,\infty]$ is narrowly lower semicontinuous.
The interior dissipation potential $\calR_\Omega: \calM_{\geq0}(\Omega)\times\calM(\Omega;\R^d)\to[0,\infty]$ is jointly lower
semicontinuous w.r.t.\ the narrow topology in both arguments.
\end{lemma}
\begin{proof}
The lower semicontinuity of the Fisher information~\eqref{eq:def:fisher} can be observed from its dual representation as follows.
From the bound $\calI(\rho|\pi)<\infty$, we obtain that $u= \dx \rho / \dx \pi \in W^{1,1}(\Omega,\pi)$ by Cauchy--Schwarz
\[
    \int_\Omega \abs*{\nabla u} \dx{\pi}=  \int_\Omega \abs*{\frac{\nabla u}{u}} u \dx{\pi}
    =  \int_\Omega \abs*{\frac{\nabla u}{u}} \dx{\rho}\leq \sqrt{2 \calI(\rho|\pi)} \sqrt{\rho(\Omega) } < \infty \,.
\]
With this, and recalling that $\pi$ has a smooth density, we have for any $\Phi\in C_c^1(\Omega,\R^d)$ and $u\in W^{1,1}(\Omega,\pi)$ the integration by parts identity
\[
    \int_\Omega \nabla u \cdot \Phi \dx{\pi} = \int_\Omega \bra*{- \nabla \log \pi  \cdot \Phi - \div  \Phi} \, u \dx{\pi} \,.
\]
Hence, dualizing the square function, we derive the dual formulation as
	\begin{align}
		\calI(u) &= \sup_{\Phi\in C_c^1(\Omega,\R^d)}  \int_\Omega  \bra*{ \frac{\nabla u}{u} \cdot \Phi - \frac{|\Phi|^2}{2}} u \dx \pi \nonumber \\
        &= \sup_{\Phi\in C_c^1(\Omega,\R^d)} \set[\bigg]{ \int_\Omega  \nabla u \cdot \Phi \dx\pi  -\int_\Omega \frac{|\Phi|^2}{2} u \dx \pi } \nonumber \\
        &= \sup_{\Phi\in C_c^1(\Omega,\R^d)} \int_\Omega  \bra[\bigg]{ -\nabla \log \pi \cdot \Phi - \div \Phi - \frac{|\Phi|^2}{2}} u \dx \pi  \,, \label{eq:dual:Fisher}
	\end{align}
    where we used that $\log \pi$ is differentiable.
	The dual formulation~\eqref{eq:dual:Fisher} shows that 
    $\calI(\cdot)$ is a supremum of narrowly continuous affine functionals on $\calM_{\ge 0}(\Omega)$.

    The lower semicontinuity of $\calR_\Omega$ follows directly from the joint convexity of its integrand, cf., e.g.,~\cite[Theorem 3.4.3]{Buttazzo1989}.
\end{proof}
\begin{proposition}[Lower semicontinuity of $\sfM$]\label{prop:M:lsc}
Let $\gamma_n,\gamma\in\calM_{\geq0}(\Gamma)$ with
$\gamma_n\to\gamma$ narrowly in $\calM(\Gamma)$, and let
$\jb_n,\jb\in(C^\infty(\Gamma))^\ast$ satisfy
$\skp{\jb_n,\psi}_\Gamma\to\skp{\jb,\psi}_\Gamma$ for every
$\psi\in C^\infty(\Gamma)$. Then, we have
\begin{equation}\label{eq:M:lsc}
    \sfM(\gamma,\jb) \;\leq\; \liminf_{n\to\infty}\ \sfM(\gamma_n,\jb_n) \,.
\end{equation}
\end{proposition}
\begin{proof}
W.l.o.g.\ let $L\coloneqq\liminf_{n\to\infty}\sfM(\gamma_n,\jb_n)<\infty$ and, passing to
a subsequence, let $\sfM(\gamma_n,\jb_n)\to L$. 
Choose $(u_n,J_n)\in M(\gamma_n,\jb_n)$ with $\calR_\Omega(u_n\pi,J_n)+\calI(u_n)\leq
\sfM(\gamma_n,\jb_n)+\tfrac1n\leq L+2$. 
Cauchy--Schwarz and the uniform bound on $\calI$ imply that the sequence $(u_n)_n$ is bounded in $W^{1,1}(\Omega,\pi)$ and $(\sqrt{u_n})_n$ is bounded in $H^1(\Omega,\pi)$.
Furthermore, bounding by $\calR_\Omega$, the sequence $(J_n)_n$ is bounded in $\calM(\Omega;\R^d)$. 
Thus, along another subsequence, $u_n\rightharpoonup u$ weakly in $W^{1,1}(\Omega,\pi)$ and therefore $\tr u_n\to\tr u$ narrowly in $\calM(\Gamma)$.
The bounds $c\leq u_n\leq C$ persist, so $u\in\Acc$, and, by linearity of~\eqref{eq:u-gamma_coupling}, it is also satisfied in the limit.

For the flux, we have $J_n\to J$ narrowly in $\calM(\Omega;\R^d)$ and $J\ll u \pi$ (again due to the uniform bound on $\calR_\Omega$ and its lower semicontinuity from Lemma~\ref{lem:FI:lower:semicont}), the linear constraint $\div J=0$ passes to the limit by testing with $C_c^\infty(\Omega)$, and for $\varphi\in C^\infty(\overline\Omega)$ we have
\[
    \int_\Omega\nabla\varphi\cdot\dx J = \lim_{n\to \infty}\int_\Omega\nabla\varphi\cdot\dx J_n = \lim_{n\to \infty}\skp{\jb_n,\tr\varphi}_\Gamma = \skp{\jb,\tr\varphi}_\Gamma \,.
\]
Hence, in total $(u,J)\in M(\gamma,\jb)$. 
This concludes the proof, since by Lemma~\ref{lem:FI:lower:semicont} we have
$\sfM(\gamma,\jb)\leq\calR_\Omega(u\pi,J)+\calI(u)\leq L$.
\end{proof}

\subsection{Compactness}\label{sec:comp}

Throughout this section we consider families $(\mu^\eps,\gamma^\eps,J^\eps,j^\eps)_{\eps>0}\subset\CE_T$
satisfying the natural a priori bounds
\begin{equation}\label{eq:uniform_bounds}
    \sup_{\eps>0}\ \calD^\eps_T(\mu^\eps,\gamma^\eps;J^\eps,j^\eps) \eqqcolon C_0 <\infty  
    \quad\text{and}\quad
    \sup_{\eps>0}\ \esssup_{t\in[0,T]}\ \calE^\eps(\mu^\eps_t,\gamma^\eps_t) \eqqcolon C_1 <\infty
    \,.
\end{equation}
In addition, we make the standing nondegeneracy assumption
\begin{equation}\label{hyp:P}
    u^\eps_t \coloneqq \frac{\dx\rho^\eps_t}{\dx\pi}\in[c,C]
    \quad\text{and}\quad
    \ub^\eps_t \coloneqq \frac{\dx\gamma^\eps_t}{\dx\omega}\in[c,C]
    \quad\text{for a.e.\ }t\in(0,T),\ \eps>0,
\end{equation}
where we recall $\rho^\eps=\mu^\eps/\eps$.
Assumption~\eqref{hyp:P} is strictly stronger than the uniform energy bound in~\eqref{eq:uniform_bounds}. 
By Lemma~\ref{lem:comparison} the system~\eqref{eq:system:eps} propagates~\eqref{hyp:P}, so that along
solutions it is a condition on the initial data alone.
Since we do not address the well-posedness of~\eqref{eq:system:eps} in this note, we keep it as a hypothesis on the family.
Nevertheless, we present proofs with the aim of dropping~\eqref{hyp:P} in the future in mind.

\begin{theorem}[Result~\ref{result:EDPlimit} --- compactness]\label{thm:comp}
Let $(\mu^\eps,\gamma^\eps,J^\eps,j^\eps)_{\eps>0}\subset \CE_T$
satisfy~\eqref{eq:uniform_bounds} and~\eqref{hyp:P}. 
Then, there exist a density
$u\in L^1(0,T;W^{1,1}(\Omega,\pi))$ with $u_t\in\Acc$ for a.e.~$t$, a measure
$J\in\calM([0,T]\times\Omega;\R^d)$ with $J\ll u(\scrL^1\otimes\pi)$,
and a pair $(\gamma,\jb)\in\lCE_T$ with
$\gamma_t = \ub_t\omega$ for $\ub_t\in\Acc_\Gamma$, such that 
$(u_t,J_t)\in M(\gamma_t,\jb_t)$ for a.e.\ $t\in[0,T]$ and
along a subsequence
\begin{enumerate}[label=(\roman*)]
    \item\label{prop:comp:mu} $\norm{\mu^\eps}_{\mathrm{TV}([0,T]\times\Omega)} \to 0$;
    \item\label{prop:comp:u} $u^\eps\rightharpoonup u$ weakly in $L^1(0,T;W^{1,1}(\Omega,\pi))$,
    and $\rho^\eps\to \rho\coloneqq u\,(\scrL^1\otimes\pi)$ narrowly in $\calM([0,T]\times\Omega)$;
    \item\label{prop:comp:gamma} $\ub^\eps \rightharpoonup \ub$ weakly in $L^1((0,T)\times\Gamma)$, with $\ub_t = \tr u_t$ for a.e.\ $t$;
    \item\label{prop:comp:J} $J^\eps \to J$ narrowly in $\calM([0,T]\times\Omega;\R^d)$.
\end{enumerate}
\end{theorem}

\begin{proof}
Since all four terms of $\calD^\eps_T$ in~\eqref{eqdef:Deps} are non-negative,
\eqref{eq:uniform_bounds} controls each of them separately, i.e., we have
\begin{equation}\label{eq:comp:interior:bounds}
    \int_0^T \calI(u^\eps_t)\dx t \leq C_0 \,,
    \qquad
    \int_0^T\!\!\int_\Omega \frac12\abs*{\frac{\dx J^\eps_t}{\dx\rho^\eps_t}}^2\dx\rho^\eps_t\dx t \leq C_0 \,,
\end{equation}
and the boundary slope gives
\begin{equation}\label{eq:comp:hellinger}
    \int_0^T\!\!\int_\Gamma \abs[\big]{\sqrt{\tr u^\eps_t} - \sqrt{\ub^\eps_t}}^2 \dx\sqrt{\tr\pi\,\omega}\dx t
    \leq \frac{\sqrt\eps}{2}\,C_0 \overset{\ \eps\to0}{\longrightarrow} 0 \,.
\end{equation}

\smallskip\noindent
\ref{prop:comp:mu}: 
We denote $\widetilde u^\eps \coloneqq \tfrac{\dx\mu^\eps}{\dx(\scrL^1\otimes\pi)}$ and recall that the domain part of the energy satisfies the bound
\begin{align*}
    \sup_{\eps>0}\esssup_{t\in[0,T]} \int_\Omega \eps \lambda_B\bra[\Big]{\frac{1}{\eps}\widetilde u^\eps(t)}\dx\pi <\infty,
\end{align*}
where $\lambda_B(r)\coloneqq r\log r -r +1$ grows superlinearly.
In particular, the family $(\widetilde u^\eps)_\eps$ is uniformly integrable and converges pointwise a.e. on $[0,T]\times\Omega$ to zero. Hence, by the Vitali convergence theorem, $\widetilde u^\eps \to 0$ in $L^1([0,T]\times\Omega)$ and thus the claim.

\ref{prop:comp:u} and \ref{prop:comp:J}: Using the $\eps$-uniform control over $\nabla \sqrt{u^\eps}$ and $J^\eps/\sqrt{u^\eps}$ in $L^2((0,T);L^2(\pi))$, this follows by arguing as in \cite[Theorem~3.2, proof step 6]{ArnrichMielkePeletierSavareVeneroni2012}.
Here, we require~\eqref{hyp:P} only to ensure absolute continuity in time w.r.t.\ the Lebesgue measure of the limit curves $\nabla u(\scrL^1\otimes\pi)$ and $J$,
and to obtain $u_t\in\Acc$.

\ref{prop:comp:gamma}: The weak convergence follows directly from the uniform energy control and the Dunford--Pettis theorem. 
The membership $\ub\in\Acc_\Gamma$ then follows from the weak convergence and~\eqref{hyp:P}.
The trace identity follows since $\tr u^\eps\rightharpoonup\tr u$ weakly in $L^1((0,T)\times\Gamma)$ by~\ref{prop:comp:u} and the continuity of the trace, so that~\eqref{eq:comp:hellinger} and the narrow lower semicontinuity of the Hellinger distance give $\tr u=\ub$.

It remains to show that the limit continuity equation holds.
To this end, testing~\eqref{eq:CE:weak} with $\varphi\in C^1_c((0,T)\times\overline\Omega)$ and
$\phi=\tr\varphi$ removes the interface term, thus yielding~\eqref{eq:CE:embedded}.
Passing to the limit with \ref{prop:comp:mu}, \ref{prop:comp:gamma}
and~\ref{prop:comp:J} we obtain~\eqref{eq:limit-CE:J}, which due to the absolute continuity in time w.r.t.\ the Lebesgue measure of $J$ also yields the narrow absolute continuity of $\gamma$.
Finally, choosing in~\eqref{eq:limit-CE:J} $\varphi\in
C_c^\infty((0,T)\times\Omega)$, we find $\div J_t=0$ for a.e.\ $t$, so that~\eqref{eq:Gauss-Green} provides the normal traces $\jb_t$ satisfying
$(\gamma,\jb)\in\lCE_T$, where we recall $\gamma_t=\ub_t\omega$. 
With these traces, we also have $(u_t,J_t)\in M(\gamma_t,\jb_t)$. 
\end{proof}

\subsection{The lower limit inequality}\label{sec:liminf}

Having established convergence in Theorem~\ref{thm:comp}, we now turn our focus to the lower limit inequality connecting the prelimit gradient systems to the limit evolution of cell problems.
We denote by $\sfE(\gamma) \coloneqq \calH(\gamma\mid\omega)$ the limit energy, and for
$(\gamma,\jb)\in\lCE_T$ we define the limit dissipation functional
\begin{equation}\label{eqdef:D0}
    \calD^0_T(\gamma,\jb) \coloneqq \int_0^T \sfM(\gamma_t,\jb_t)\dx t \,.
\end{equation}
The integral in~\eqref{eqdef:D0} is well defined, since the integrand is Lebesgue measurable.
Indeed, thanks to $t\mapsto\abs{J_t}(\Omega)$ being integrable, 
the sets
$A_N\coloneqq\set{t\in[0,T] : \abs{J_t}(\Omega)\leq N}$ exhaust $[0,T]$ up to a null set.
For $t\in A_N$, the narrow continuity
of $\gamma$ and the bound~\eqref{eq:jb:bound} show that the pair $(\gamma_t,\jb_t)$ lives in
a compact metrizable subset of $\calM_{\geq0}(\Gamma)\times\bigl(C^\infty(\Gamma)\bigr)^\ast$.
On this set $\sfM$ is lower semicontinuous by Proposition~\ref{prop:M:lsc}, hence Borel.
The map $t\mapsto(\gamma_t,\jb_t)$ is measurable, since $\gamma$ is continuous and the
pairing $t\mapsto\skp{\jb_t,\tr\varphi}_\Gamma=\int_\Omega\nabla\varphi\cdot\dx J_t$ is measurable. 
Hence $t\mapsto\sfM(\gamma_t,\jb_t)$ is measurable on every $A_N$, and thus on $[0,T]$.

\begin{theorem}[Result~\ref{result:EDPlimit} --- lower limit inequality]\label{thm:liminf}
Let $(\mu^\eps,\gamma^\eps, J^\eps, j^\eps)_{\eps>0}\subset \CE_T$
satisfy~\eqref{eq:uniform_bounds} and~\eqref{hyp:P}, and let
$(\gamma,\jb)\in\lCE_T$ be its limit as in Theorem~\ref{thm:comp}. 
Then, 
we have 
the lower limit inequalities
\begin{align}\label{eq:liminf:diss}
    \liminf_{\eps\to 0} \calD_T^\eps(\mu^\eps,\gamma^\eps,J^\eps,j^\eps) &\geq \calD_T^0(\gamma,\jb)\,,\\
    \label{eq:liminf:energy}
    \liminf_{\eps\to 0} \calE^\eps(\mu^\eps_t,\gamma^\eps_t) &\ge \sfE(\gamma_t) \quad\text{for all $t\in[0,T]$}.
\end{align}
In particular, if the initial data are \emph{well prepared}, i.e.,
\begin{equation}\label{hyp:liminf:wellprep}
    \calE^\eps(\mu^\eps_0,\gamma^\eps_0) \longrightarrow \sfE(\gamma_0)
    \qquad\text{as }\eps\to0\,,
\end{equation}
then
\begin{equation}\label{eq:liminf:L}
    \liminf_{\eps\to 0}\calL^\eps_T(\mu^\eps,\gamma^\eps;J^\eps,j^\eps) \geq \calL^{\DtN}_T(\gamma)
\end{equation}
with $\calL^{\DtN}_T$ from~\eqref{eq:DeGiorgi:red}. 
\end{theorem}
\begin{proof}
Dropping the non-negative boundary terms on the prelimit side, we have
\[
    \calD^\eps_T(\mu^\eps,\gamma^\eps;J^\eps,j^\eps)
    \geq \int_0^T \bigl[ \calR_\Omega(\rho^\eps_t,J^\eps_t) + \calI(u^\eps_t) \bigr] \dx{t} \,.
\]
In particular, the $\eps$-dependence is only in the arguments of $\calR_\Omega$, resp.\ $\calI$.
Therefore, since $\calR_\Omega$ and $\calI$ are narrowly lower semicontinuous by Lemma~\ref{lem:FI:lower:semicont}, Fatou's lemma yields
\begin{equation}\label{eq:liminf:interior}
    \liminf_{\eps\to0} \calD^\eps_T(\mu^\eps,\gamma^\eps;J^\eps,j^\eps)
    \geq \int_0^T \bigl[ \calR_\Omega(\rho_t,J_t) + \calI(u_t) \bigr] \dx t \,.
\end{equation}
Since $(u_t,J_t)\in M(\gamma_t,\jb_t)$ for a.e.\ $t$ by Theorem~\ref{thm:comp}, we can further bound from below by taking the infimum over $M(\gamma_t,\jb_t)$ pointwise in $t$
to obtain~\eqref{eq:liminf:diss}.

For the energy, we drop the term $\calH(\mu\mid\pi^\eps)\geq0$, obtaining
$\calE^\eps(\mu^\eps_t,\gamma^\eps_t)\geq\calH(\gamma^\eps_t\mid\omega)$.
Then, narrow lower semicontinuity of the relative entropy yields for all $t\in[0,T]$
\begin{equation*}
    \liminf_{\eps\to0}\calE^\eps(\mu^\eps_t,\gamma^\eps_t) \geq \calH(\gamma_t\mid\omega) = \sfE(\gamma_t)\,.
\end{equation*}
Finally, combining this with the well-preparedness~\eqref{hyp:liminf:wellprep}, we obtain~\eqref{eq:liminf:L}.
\end{proof}

\section{The cell problem characterizes a gradient system}\label{sec:CELL:GS}

We now show that the cell formula~\eqref{eq:CELL} splits into a primal and a
dual dissipation potential, evaluated at the driving force of the boundary
entropy. Throughout this section $\ub\in\Acc_\Gamma$ is a fixed boundary
density, $\gamma = \ub \omega$; we recall the definition of $\sfE$,
\begin{equation}\label{eq:def:E:limit}
    \sfE(\gamma) \coloneqq \calH(\gamma\mid\omega) = \int_\Gamma \bigl(\ub\log\ub - \ub + 1\bigr)\dx\omega \,,
    \qquad \rmD\sfE(\gamma) = \log\ub \,.
\end{equation}

We begin by introducing homogeneous weighted half-order spaces, which turn out to be well suited to our problem.
For $u\in\Acc$ and $\eta\in H^{1/2}(\Gamma)$ we define the homogeneous trace seminorm
\begin{equation}\label{eq:def:H12}
    \norm{\eta}_{\dot H^{1/2}(\Gamma,u\pi)}^2
    \coloneqq \min \set*{ \int_\Omega \abs{\nabla \xi}^2 u \dx\pi \,:\, \xi \in H^1(\Omega,u\pi),\ \tr\xi = \eta } \,,
\end{equation}
the minimum being attained by the $u\pi$-harmonic extension of $\eta$. 
For $u\equiv1$ this is the classical $\pi$-weighted $H^{1/2}$-seminorm
$\norm{\cdot}_{\dot H^{1/2}(\Gamma,\pi)}$. Since $c\leq u\leq C$,
\begin{equation}\label{eq:H12:equivalence}
    c\,\norm{\eta}^2_{\dot H^{1/2}(\Gamma,\pi)} \leq \norm{\eta}^2_{\dot H^{1/2}(\Gamma,u\pi)} \leq C\,\norm{\eta}^2_{\dot H^{1/2}(\Gamma,\pi)} \,,
\end{equation}
so all these seminorms are equivalent, and they vanish exactly on the constants.
Accordingly we let
\begin{equation}\label{eq:def:H-12}
    \begin{split}
    \dot H^{-1/2}(\Gamma) &\coloneqq \set*{\jb \in \bigl(H^{1/2}(\Gamma)\bigr)^* : \
            \skp{\jb,1}_\Gamma = 0} \,,\\
    \tfrac{1}{2}\norm{\jb}^2_{\dot H^{-1/2}(\Gamma,u\pi)} &\coloneqq 
            \sup_{\eta\in H^{1/2}(\Gamma)}
            \bigl\{ 
                \skp{\jb,\eta}_\Gamma 
                - \tfrac{1}{2} \norm{\eta}_{\dot H^{1/2}(\Gamma,u\pi)}^2
            \bigr\}
            \,,
    \end{split}
\end{equation}
the mean-zero restriction removing the kernel of the seminorm~\eqref{eq:def:H12}.
By~\eqref{eq:H12:equivalence} the norms in~\eqref{eq:def:H-12} are again all
equivalent, so that finiteness of $\norm{\jb}_{\dot H^{-1/2}(\Gamma,u\pi)}$ does not
depend on $u$.

We observe that among all extensions of $\log\ub$, the logarithm of the harmonic extension is the one of least weighted energy, and its energy is twice the Fisher information.

\begin{lemma}[Optimality of the harmonic extension]\label{lem:logu:optimal}
Let $\ub\in\Acc_\Gamma$ and $u \coloneqq \calU(\ub)$. Then $u\in\Acc$ and $\log u\in
H^1(\Omega,u\pi)$ with $\tr\log u = \log\ub$, and $\log u$ is the unique
minimizer of the $u\pi$-weighted Dirichlet energy among the extensions of
$\log\ub$, that is
\begin{equation}\label{eq:logu:optimal}
    \int_\Omega\abs{\nabla\log u}^2 u \dx\pi \leq \int_\Omega\abs{\nabla\xi}^2 u \dx\pi
    \qquad\text{for all } \xi\in H^1(\Omega,u\pi) \text{ with } \tr\xi = \log\ub \,,
\end{equation}
with equality if and only if $\xi = \log u$. In particular
\begin{equation}\label{eq:FisherOptimality}
    \norm{\log\ub}^2_{\dot H^{1/2}(\Gamma,\calU(\ub)\pi)}
    = \int_\Omega\abs{\nabla\log \calU(\ub)}^2 \calU(\ub)\dx\pi
    = 2\,\calI\bigl(\calU(\ub)\bigr) \,.
\end{equation}
\end{lemma}
\begin{proof}
By the maximum principle, we immediately have $u=\calU(\ub)\in\Acc$.
Hence $c\leq u\leq C$ and so $s\mapsto\log s$ is
Lipschitz on $[c,C]$, which implies  $\log u\in H^1(\Omega,u\pi)$ with $\nabla\log u =
\nabla u/u$ and $\tr\log u = \log\tr u = \log\ub$.
In particular
\begin{equation}\label{eq:fisher-id}
    \int_\Omega\abs{\nabla\log u}^2 u\dx\pi = \int_\Omega\frac{\abs{\nabla u}^2}{u}\dx\pi = 2\calI(u) \,.
\end{equation}
On the domain $\Omega$ the trace-free elements of $H^1(\Omega,u\pi)$ are exactly $H^1_0(\Omega,u\pi) = H^1_0(\Omega)$, since $u \pi$ is equivalent to the Lebesgue measure. 
So every competitor in~\eqref{eq:logu:optimal} is of the form $\xi = \log u+\zeta$ with $\zeta\in H^1_0(\Omega,u\pi)$. 
Using $u\nabla\log u = \nabla u$, we expand and estimate
\begin{align*}
    \int_\Omega\abs{\nabla\log u + \nabla\zeta}^2 u\dx\pi
    &= \int_\Omega\abs{\nabla\log u}^2 u\dx\pi + \int_\Omega\abs{\nabla\zeta}^2u\dx\pi \\
    &\qquad + 2\int_\Omega\nabla u\cdot\nabla\zeta\dx\pi
    \;\geq\; \int_\Omega\abs{\nabla\log u}^2 u\dx\pi \,,
\end{align*}
the cross term vanishing because it is the weak formulation of
$\div(\pi\nabla u)=0$ tested with $\zeta\in H^1_0(\Omega,\pi)$. Equality forces
$\int_\Omega\abs{\nabla\zeta}^2u\dx\pi = 0$, hence $\nabla\zeta=0$ by $u\geq c$,
so $\zeta$ is constant and therefore, being trace-free, $\zeta=0$. This
is~\eqref{eq:logu:optimal}, and together with~\eqref{eq:fisher-id} and the
definition~\eqref{eq:def:H12} of the trace norm also~\eqref{eq:FisherOptimality}.
\end{proof}

Our next step is to derive from the cell formula a primal-dual pair of dissipation potentials. This is what we do in the following two lemmas.  
\begin{lemma}[Properties of the cell formula]\label{lem:props:cell}
    Let $\ub\in\Acc_\Gamma$ and $\gamma = \ub\,\omega$. 
    Then, the no-flux cell
    problem has the value
    \begin{equation}\label{eq:char:M-zero-dot-gamma}
        \sfM(\gamma,0) = 2\norm{\sqrt{\ub}}^2_{\dot H^{1/2}(\Gamma,\pi)} <\infty \,,
    \end{equation}
    and the functional
	\begin{equation}\label{eq:primal:cell}
	\sfR(\gamma, \jb) \coloneqq \sfM(\gamma,\jb) - \sfM(\gamma,0)\,,
    \qquad \jb\in\bigl(C^\infty(\Gamma)\bigr)^*,
	\end{equation}
	is well defined and an \emph{admissible dissipation functional}, that is
	\begin{enumerate}[label=(\roman*)]
		\item $\sfR(\gamma,\jb)\geq 0$ for all $\jb$, and $\sfR(\gamma,0)=0$;
		\item $\jb \mapsto \sfR(\gamma,\jb)$ is convex.
	\end{enumerate}
	Moreover, the cell formula
	$(\gamma,\jb)\mapsto\sfM(\gamma,\jb)$ is jointly convex on
	$\calM_{\geq0}(\Gamma)\times\bigl(C^\infty(\Gamma)\bigr)^*$.
\end{lemma}
\begin{proof}
For $\jb=0$ the choice $J\equiv0$ is admissible in~\ref{def:cell:i} and
optimal, since $\calR_\Omega\geq0$ and $\calR_\Omega(u\pi,0)=0$. Hence
$\sfM(\gamma,0) = \inf\set{\calI(u) : u\in\Acc,\ \tr u=\ub} 
$.

Substituting
 $\xi\coloneqq\sqrt u$, which maps $\set{u\in\Acc:\tr u=\ub}$ bijectively onto
 $\set{\xi\in H^1(\Omega,\pi):\sqrt c\leq\xi\leq\sqrt C,\ \tr\xi=\sqrt\ub}$, this
 infimum equals $2\inf\int_\Omega\abs{\nabla\xi}^2\dx\pi$ over that set. The
 minimizer of the unconstrained problem is $\calU(\sqrt\ub)$, which satisfies
 $\sqrt c\leq\calU(\sqrt\ub)\leq\sqrt C$ by the maximum principle applied to the datum
 $\sqrt\ub\in[\sqrt c,\sqrt C]$.
 Hence, the constraint is inactive
 and~\eqref{eq:char:M-zero-dot-gamma} just becomes~\eqref{eq:def:H12}. 
 Since $\ub\in \Acc_\Gamma$ and the square-root is Lipschitz on $[c,C]$, we have $\sqrt\ub\in H^{1/2}(\Gamma)$, i.e., we also have the finiteness in~\eqref{eq:char:M-zero-dot-gamma}.

Next, dropping the term $\calR_\Omega(u\pi,J)\geq0$ in~\eqref{eq:CELL} gives
$\sfM(\gamma,\jb)\geq\sfM(\gamma,0)$ for every $\jb$, i.e., $\sfR\geq0$ and
$\sfR(\gamma,0)=0$ by definition. 

For the convexity, let $\jb_1,\jb_2$ and
$\lambda_1,\lambda_2\in(0,1)$ with $\lambda_1+\lambda_2=1$. If $(u_i,J_i)\in
M(\gamma,\jb_i)$, then $(\lambda_1u_1+\lambda_2u_2,\lambda_1J_1+\lambda_2J_2)
\in M(\gamma,\lambda_1\jb_1+\lambda_2\jb_2)$.
Indeed, the set $\Acc$ is convex, the convex combination $\lambda_1J_1+\lambda_2J_2$ is still divergence-free, and by the linearity of the trace, the trace condition~\eqref{eq:u-gamma_coupling} is preserved.
From this, we conclude by observing that the integrand
in~\eqref{eq:CELL} is jointly convex in $(u,J)$, 
i.e., $\jb\mapsto\sfM(\gamma,\jb)$ is convex and hence so is $\sfR(\gamma,\cdot)$.

For the joint convexity, let $\gamma_i=\ub_i\omega$ for $i=1,2$.
If $\sfM(\gamma_i,\jb_i)=+\infty$ for some~$i$, there is nothing to show.
Otherwise, the same argument applies to $(u_i,J_i)\in M(\gamma_i,\jb_i)$.
Indeed, the trace condition~\eqref{eq:u-gamma_coupling} is linear
in the pair $(u,\ub)$,
so that $\lambda_1u_1+\lambda_2u_2$ has the trace
$\lambda_1\ub_1+\lambda_2\ub_2$.
\end{proof}

Having defined the primal dissipation potential $\sfR$ in the previous lemma, we define the dual dissipation potential as its convex conjugate w.r.t.\ the second argument,
\begin{equation}\label{eq:def:Rstar}
    \sfR^*(\gamma,\eta) \coloneqq \sup\set[\big]{\skp{\jb,\eta}_\Gamma - \sfR(\gamma,\jb) \;:\; \jb\in \dot H^{-1/2}(\Gamma)} \,.
\end{equation}
We now further characterize the primal and dual dissipation potentials. 
\begin{lemma}[Characterizations of $\sfM$, $\sfR$ and $\sfR^*$]\label{lem:M:R:Rstar:char}
Let $\ub\in\Acc_\Gamma$, $\gamma = \ub\omega$, and $\jb\in \dot H^{-1/2}(\Gamma)$.
\begin{enumerate}[label=(\roman*)]
\item\label{lem:MRR:i} For every $u\in\Acc$ with $\tr u=\ub$,
\begin{equation}\label{eq:min-pb-over-J-equals-norm-of-jb}
    \min \set[\bigg]{ \calR_\Omega(u\pi,J) \,\bigg|\, J\in \calM(\Omega;\R^d): \  \div J = 0, \ J\cdot n = \jb }
    = \frac12 \norm{\jb}^2_{\dot H^{-1/2}(\Gamma,u\pi)} \,.
\end{equation}
\item\label{lem:MRR:ii} With all infima and suprema taken over $u\in\Acc$ with $\tr u=\ub$,
\begin{align}
        \sfM(\gamma,\jb) &= \inf_{u: \tr u=\ub} \set*{ \frac{1}{2} \norm{\jb}_{\dot H^{-1/2}(\Gamma,u \pi)}^2 + \calI(u) } \,, \label{eq:char-M-jb-HMOH}\\
        \sfR(\gamma,\jb) &= \inf_{u: \tr u=\ub} \set*{ \frac{1}{2} \norm{\jb}_{\dot H^{-1/2}(\Gamma,u \pi)}^2 + \calI(u) - 2\norm{\sqrt{\ub}}^2_{\dot H^{1/2}(\Gamma,\pi)} } \,, \label{eq:char-R-jb-HMOH}\\
        \sfR^*(\gamma,\eta) &= \sup_{u: \tr u =\ub} \set*{ \frac{1}{2} \norm{\eta}_{\dot H^{1/2}(\Gamma,u \pi)}^2 - \calI(u) + 2\norm{\sqrt{\ub}}^2_{\dot H^{1/2}(\Gamma,\pi)} } \,,\label{eq:char-Rstar-jb-HMOH}
    \end{align}
the last identity for every $\eta\in H^{1/2}(\Gamma)$.
\end{enumerate}
\end{lemma}
\begin{proof}
\ref{lem:MRR:i}: For every admissible $J$ 
and every $\phi\in C^\infty_c(\R^d)$, the Cauchy--Schwarz inequality gives
\begin{equation}\label{eq:MRR:CS}
    \abs[\big]{\skp{\jb,\tr\phi}_\Gamma}
    = \abs[\bigg]{\int_\Omega \nabla \phi \cdot \dx J}
    \leq \bra[\bigg]{\int_\Omega \abs{\nabla \phi}^2 \,u\dx \pi}^{1/2}
    \bigl(2\calR_\Omega(u\pi,J)\bigr)^{1/2} ,
\end{equation}
so by~\eqref{eq:def:H12} and~\eqref{eq:def:H-12} the left-hand side
of~\eqref{eq:min-pb-over-J-equals-norm-of-jb} is at least as large as the right-hand side.
For the converse, let $\overline\calV$ be the closure in $L^2(\Omega,u\pi)$ of
$\calV\coloneqq \set[\big]{\nabla \phi|_{\overline \Omega}: \phi\in C^\infty_c (\R^d)}$
and consider $T:\calV\to\R$, $T\nabla\phi \coloneqq \skp{\jb,\tr\phi}_\Gamma$, which is
well defined because $\skp{\jb,1}_\Gamma=0$, so that $T$ vanishes on constants.
By~\eqref{eq:def:H12} and~\eqref{eq:def:H-12},
\begin{align}\label{eq:proof:lem:M:R:Rstar:char:T:bounded}
    \abs{T\nabla\phi} \leq \norm{\jb}_{\dot H^{-1/2}(\Gamma,u\pi)}\norm{\tr\phi}_{\dot H^{1/2}(\Gamma,u\pi)}
    \leq \norm{\jb}_{\dot H^{-1/2}(\Gamma,u\pi)}\Bigl(\int_\Omega\abs{\nabla\phi}^2u\dx\pi\Bigr)^{1/2} ,
\end{align}
so $T$ is a bounded linear functional on the Hilbert space
$(\overline\calV,\norm{\cdot}_{L^2(\Omega,u\pi)})$ and its norm is exactly
$\norm{\jb}_{\dot H^{-1/2}(\Gamma,u\pi)}$, since the traces of $C_c^\infty(\R^d)$-functions
are dense in $H^{1/2}(\Gamma)$. 
Hence, by Riesz' representation theorem, there exists $f\in\overline\calV$ with
$\norm{f}_{L^2(\Omega,u\pi)} = \norm{\jb}_{\dot H^{-1/2}(\Gamma,u\pi)}$ and
$\skp{\jb,\tr\phi}_\Gamma = \int_\Omega f\cdot \nabla \phi \, u\dx \pi$ for all
$\phi\in C_c^\infty(\R^d)$.
We set $J \coloneqq fu\pi$. It is a finite measure, since $\int_\Omega\abs{f}u\dx\pi\leq\norm{f}_{L^2(\Omega,u\pi)}\,(u\pi(\Omega))^{1/2}<\infty$.
By the representation identity above, it satisfies $\div J=0$ with $J\cdot n=\jb$.
Finally, $\calR_\Omega(u\pi,J) = \tfrac12\norm{f}^2_{L^2(\Omega,u\pi)}$, so the minimum is
attained and the two sides agree. 

\ref{lem:MRR:ii}: Splitting the infimum in~\eqref{eq:CELL} as $\inf_u\inf_J$ and
inserting~\eqref{eq:min-pb-over-J-equals-norm-of-jb} for the inner one
gives~\eqref{eq:char-M-jb-HMOH}. Subtracting~\eqref{eq:char:M-zero-dot-gamma},
which is finite, yields~\eqref{eq:char-R-jb-HMOH}. For~\eqref{eq:char-Rstar-jb-HMOH}, by definition
\begin{align*}
&\sfR^*(\gamma,\eta) 
= \sup_{\substack{u: \tr u =\ub\\ \jb\in \dot H^{-1/2}(\Gamma)}}
  \set*{\skp {\jb,\eta}_\Gamma - \frac{1}{2} \norm{\jb}_{\dot H^{-1/2}(\Gamma,u \pi)}^2 - \calI(u) + 2\norm{\sqrt{\ub}}^2_{\dot H^{1/2}(\Gamma,\pi)}}\\
&= \sup_{u: \tr u =\ub} \set[\Bigg]{\ \sup_{\jb\in \dot H^{-1/2}(\Gamma)} \set*{\skp {\jb,\eta}_\Gamma - \frac{1}{2} \norm{\jb}_{\dot H^{-1/2}(\Gamma,u \pi)}^2} - \calI(u) + 2\norm{\sqrt{\ub}}^2_{\dot H^{1/2}(\Gamma,\pi)}}\,.
\end{align*} 
The inner supremum is
$\tfrac12\norm{\eta}^2_{\dot H^{1/2}(\Gamma,u\pi)}$, 
by Legendre duality of the seminorm in~\eqref{eq:def:H-12}.
\end{proof}

\begin{remark}[Quadratic bounds]\label{rem:R:bounds}
Let $\ub\in\Acc_\Gamma$, $\gamma=\ub\omega$, and
$u_*\coloneqq\calU(\sqrt\ub)^2$.
Then, for all $\jb\in\dot H^{-1/2}(\Gamma)$ and $\eta\in H^{1/2}(\Gamma)$,
\begin{align*}
    \frac{1}{2C}\norm{\jb}^2_{\dot H^{-1/2}(\Gamma,\pi)}
    &\leq \sfR(\gamma,\jb)
    \leq \frac12\norm{\jb}^2_{\dot H^{-1/2}(\Gamma,u_*\pi)}
    \leq \frac{1}{2c}\norm{\jb}^2_{\dot H^{-1/2}(\Gamma,\pi)} \,,\\
    \frac{c}{2}\norm{\eta}^2_{\dot H^{1/2}(\Gamma,\pi)}
    &\leq \frac12\norm{\eta}^2_{\dot H^{1/2}(\Gamma,u_*\pi)}
    \leq \sfR^*(\gamma,\eta)
    \leq \frac{C}{2}\norm{\eta}^2_{\dot H^{1/2}(\Gamma,\pi)} \,.
\end{align*}
Indeed, by the proof of Lemma~\ref{lem:props:cell},
$u_*\in\Acc$ satisfies $\tr u_*=\ub$ and $\calI(u_*)=\sfM(\gamma,0)$,
while $\calI(u)\geq\sfM(\gamma,0)$ for every $u\in\Acc$ with $\tr u=\ub$.
Choosing $u=u_*$ in~\eqref{eq:char-R-jb-HMOH}
and~\eqref{eq:char-Rstar-jb-HMOH} gives the upper bound for $\sfR$
and the lower bound for $\sfR^*$.
Dropping the term $\calI(u)-\sfM(\gamma,0)\geq0$ gives the lower bound
for $\sfR$ and the upper bound for $\sfR^*$,
once we use~\eqref{eq:H12:equivalence} and its dual version
$C^{-1}\norm{\jb}^2_{\dot H^{-1/2}(\Gamma,\pi)}
\leq\norm{\jb}^2_{\dot H^{-1/2}(\Gamma,u\pi)}
\leq c^{-1}\norm{\jb}^2_{\dot H^{-1/2}(\Gamma,\pi)}$.
The outer inequalities follow directly from these two norm equivalences.
In particular, $\sfR(\gamma,\cdot)$ is finite on all of $\dot H^{-1/2}(\Gamma)$,
and the quadratic potential with the weight of the Fisher-harmonic
extension~$u_*$ bounds $\sfR$ from above and $\sfR^*$ from below.
\end{remark}

Finally, we are in the position to prove the duality splitting for the cell formula.
\begin{theorem}[Duality splitting of the cell formula --- Result~\ref{result:CELL:splitting}]\label{thm:splitting}
Let $\ub\in\Acc_\Gamma$, $\gamma = \ub \omega$,  and let $\jb\in \bigl(C^\infty(\Gamma)\bigr)^*$. Then
\begin{equation}\label{eq:splitting}
    \sfM(\gamma,\jb) = \sfR(\gamma,\jb) + \sfR^*\bigl(\gamma,-\rmD\sfE(\gamma)\bigr)\,,
\end{equation}
with $\sfE$ from~\eqref{eq:def:E:limit} and $\rmD\sfE(\gamma) = \log\ub$.
\end{theorem}
\begin{proof}
For the proof it is sufficient to show that
$\sfR^*(\gamma,-\log\ub) = \sfM(\gamma,0)$.
Indeed, the
splitting~\eqref{eq:splitting} then follows from the
definition~\eqref{eq:primal:cell} of $\sfR$.
This holds for any~$\jb$, since if $\sfM(\gamma,\jb)=+\infty$ both sides of~\eqref{eq:splitting}
are $+\infty$.
Here we use that the subtracted term $\sfM(\gamma,0)$ in~\eqref{eq:primal:cell} is finite by~\eqref{eq:char:M-zero-dot-gamma}.
Since $\ub\in\Acc_\Gamma$ and
$s\mapsto\log s$ is Lipschitz on $[c,C]$, we have $\log\ub\in H^{1/2}(\Gamma)$.
Hence, the choice $\eta=-\log\ub$ is admissible in~\eqref{eq:char-Rstar-jb-HMOH}.

We first consider the lower bound. 
Since, by Lemma~\ref{lem:logu:optimal}, $u = \calU(\ub)$ is admissible, we may estimate
\begin{equation}\label{eq:Rstar-at-log}
    \begin{split}
    \sfR^*\bigl(\gamma,-\log\ub\bigr) 
    &\geq \tfrac{1}{2}\norm{\log\ub}^2_{\dot H^{1/2}(\Gamma,\calU(\ub)\pi)} - \calI(\calU(\ub)) + 2\norm{\sqrt\ub}^2_{\dot H^{1/2}(\Gamma,\pi)}\\
    &= 2\norm{\sqrt\ub}^2_{\dot H^{1/2}(\Gamma,\pi)} \,,
    \end{split}
\end{equation}
where in the second step we used~\eqref{eq:FisherOptimality}.

For the upper bound, let $u\in\Acc$ with $\tr u=\ub$ be an arbitrary competitor.
Since $c\leq u\leq C$, we have $\log u\in H^1(\Omega,u\pi)$ with $\tr\log u =
\log\ub$, so $\log u$ is admissible in~\eqref{eq:def:H12} and it holds
\begin{equation}\label{eq:logub-vs-Fisher}
    \norm{\log\ub}^2_{\dot H^{1/2}(\Gamma,u\pi)} \leq \int_\Omega\abs{\nabla\log u}^2 u\dx\pi = 2\calI(u)\,.
\end{equation}
Consequently the expression inside the supremum
in~\eqref{eq:char-Rstar-jb-HMOH} is bounded by
$2\norm{\sqrt\ub}^2_{\dot H^{1/2}(\Gamma,\pi)}$, uniformly in $u$, which concludes the proof.
\end{proof}

\section{Discussion}\label{sec:discussion}

\subsection*{Strand 1: Cell-formula connection to one-dimensional case}

The cell formula~\eqref{eq:CELL} contains the cell problem
of~\cite{ArnrichMielkePeletierSavareVeneroni2012} as its one-dimensional case.
Take $d=1$, $\Omega=(0,1)$, $\pi=\scrL^1$, take for $\omega$ the counting measure on
$\Gamma=\set{0,1}$, and write $\ub=(\ub_0,\ub_1)$ and $\gamma = \ub_0\delta_0 + \ub_1\delta_1$. The constraint $\div J=0$ in $\Omega$
forces $J$ to be a constant $w\in\R$, so that $\jb=w(\delta_1-\delta_0)$ and the
minimization over $J$ in Definition~\ref{def:cell} is trivial. 
What remains is
\begin{equation}\label{eq:cell:1d}
    \sfM(\gamma,\jb)
    = \inf_{u\in\Acc:\ u(0)=\ub_0,\ u(1)=\ub_1}\ \frac12\int_0^1\Bigl(\frac{w^2}{u} + \frac{(u')^2}{u}\Bigr)\dx s \,,
\end{equation}
which is precisely the minimization problem defined in~\cite[Eq.~(5.1)]{ArnrichMielkePeletierSavareVeneroni2012}. 
The
no-flux value~\eqref{eq:char:M-zero-dot-gamma} evaluates here to
$\sfM(\gamma,0)=2\bigl(\sqrt{\ub_1}-\sqrt{\ub_0}\bigr)^2$ and is attained at the extension with $\sqrt u$ affine.

For two boundary points the duality splitting of $\sfM$ was identified
in~\cite{LieroMielkePeletierRenger2017} with a dual dissipation potential of
$\cosh$ type, which is what exhibits the limiting reaction as a generalized
gradient flow. Theorem~\ref{thm:splitting} provides the splitting for a general
domain, but no explicit formula for $\sfR$: in $d\geq2$ the boundary carries a
continuum of degrees of freedom and $\sfR(\gamma,\cdot)$ is the infimal-convolution-type object~\eqref{eq:def:R:intro}.
If and how this is related to an integral functional of $\cosh$ type,
we do not know.

\emph{Outlook: additivity and current fluctuations.} Writing
$\sfM(\Omega;\gamma,\jb)$ to make the domain explicit, the one-dimensional formula
obeys the series law obtained by cutting $(0,1)$ at an interior point and
optimizing over the density transmitted there. Its multidimensional counterpart
would be an \emph{additivity principle} in the sense of Bodineau and
Derrida~\cite{BodineauDerrida2004}: for a decomposition
$\overline\Omega=\overline\Omega_1\cup\overline\Omega_2$ with interface
$\Gamma_{12}$, $\sfM$ should satisfy
\[
    \sfM(\Omega;\gamma,\jb)
    = \mkern-6mu \inf_{\gamma_{12},\,\jb_{12}} \mkern-2mu  \set[\Big]{ 
        \sfM(\Omega_1;\gamma|_{\Gamma_1}\mkern-4mu +\gamma_{12},\,\jb|_{\Gamma_1}\mkern-4mu +\jb_{12})
        + \sfM(\Omega_2;\gamma|_{\Gamma_2} \mkern-4mu +\gamma_{12},\,\jb|_{\Gamma_2}\mkern-4mu  -\jb_{12})
    } .
\]
On finite
weighted graphs the corresponding series and parallel laws are proved
in~\cite[Corollary~3.4 and §7]{PeletierSchlichting2023}, see
also~\cite[§6.3]{GabrielliHarris2025}; we expect them to be recoverable from the
continuous statement by graph approximation of $\Omega$.
We expect a conceptual link
by considering the Euler--Lagrange equation for the cell problem~\eqref{eq:CELL},
which gives rise to a nonlinear
potential $\varphi[u,\jb]$ solving $\div(u\nabla\varphi[u,\jb])=0$ in $\Omega$ with
$u\,\partial_n\varphi=\jb$ on~$\Gamma$.
In terms of $\varphi[u,\jb]$ the first term in the cell formula (the flux term) becomes a
$u$-weighted capacity and $\sfM$ interpolates between the linear capacity  ($u\equiv 1$) and the Fisher information ($\jb \equiv 0$). 
Similar functionals occur as rate functionals of
boundary-driven current fluctuations in macroscopic fluctuation theory, where
additivity is a standard device for computing current
statistics~\cite{BodineauDerrida2004}. In that language,
\cite{AkkermansBodineauDerridaShpielberg2013} argues that for the symmetric
simple exclusion process connected to two reservoirs the current statistics on a
large domain in dimension~$d$ agree with those in $d=1$.
Identifying suitable boundary conditions for making the
correspondence precise is a natural next step.

\subsection*{Strand 2: Connections to diffusion processes and semigroups}

The form of the $\eps>0$ system~\eqref{eq:system:eps} of this paper suggests that it is the Fokker--Planck equation of a diffusion process in $\Omega$ combined with a reservoir on $\Gamma$. In a simplistic description, a particle diffuses in the bulk until it hits $\Gamma$, at which point it enters the reservoir with a certain probability; simultaneously, at some rate a particle that is located in the reservoir jumps into the bulk and diffuses. This heuristic description and the corresponding boundary condition~\eqref{eq:system:eps:density:Robin} have been studied in the chemical literature by e.g.\ Agmon and Szabo~\cite{AgmonSzabo90}. 

We conjecture that using capture and release based on the local time at the boundary, a process of the appropriate form can be constructed that has~\eqref{eq:system:eps} as its Fokker--Planck equation. Assuming that this is the case, this gives the Dirichlet-to-Neumann evolution a clear interpretation as the fast-diffusion limit of this process. 

Note that the setup is similar to that of `Sticky Brownian Motion', or equivalently Wentzell-boundary diffusion processes (see e.g.~\cite{Engel2003Laplacian,GrothausVoss2017,CasterasMonsaingeonSantambrogio2025,BormannMonsaingeonRengerRenesse2025}), in which a diffusing particle that hits the boundary enters the reservoir immediately, rather than with some rate. 

The identification of the $\DtN$-evolution as a scaling limit of bulk--boundary
diffusion should be compared with two neighbouring lines of work at the level of
the underlying generators. First, in~\cite{MugnoloNittka2012} it is shown that
the Wentzell--Laplace operator does \emph{not} converge to the $\DtN$-operator, not even
in the weak-resolvent sense (see~\cite[Proposition~4.15]{MugnoloNittka2012}). 
Indeed, the natural boundary-dominated rescaling of the Wentzell dynamics fails to produce
the $\DtN$-semigroup. Secondly, at the eigenvalue level,
\cite{LambertiProvenzano2017} obtains the Steklov spectrum as the singular limit
of Neumann eigenvalues under a mass concentration at the boundary, thereby
linking two classical spectral problems that are structurally close to ours. Our
$\eps\to0$ limit produces the $\DtN$-evolution directly from a bulk diffusion
coupled to a boundary reservoir.
The two references above suggest that the
intermediate Wentzell regime is genuinely different.
In our approach, the $\DtN$-generator arises in the limit from a specific scaling, namely the joint entropy
$\calE^\eps=\eps\calH(\rho^\eps|\pi)+\calH(\gamma|\omega)$ together with the
$\eps^{\pm1/2}$-weighted exchange dissipation on~$\Gamma$.

\subsection*{Strand 3: Connection to Mielke's NESS theory}

The variational convergence established by Results~\ref{result:EDPlimit}
and~\ref{result:CELL:splitting} is of EDP type, but it is \emph{not} of the
\emph{tilt-EDP} type introduced by Mielke, Montefusco and
Peletier~\cite{MielkeMontefuscoPeletier2021}. Indeed, the limit dissipation
potential depends, through the optimization in~\eqref{eq:def:R:intro}
and~\eqref{eq:def:Rs:intro} and through the integration measure $\dx\pi$ in the
weighted norms, on the reference measure~$\pi$, which in the microscopic
gradient system enters only through the driving energy. Tilt-EDP convergence
would require the limit dissipation potential to remain unchanged under
arbitrary smooth tilts of the energy, which fails here because the geometry of
the limit cell problem is dictated by the $\pi$-harmonic extension. This is
consistent with the general phenomenon observed
in~\cite{MielkeMontefuscoPeletier2021,PeletierSchlichting2023,Mielke2023,FrenzelMielke2026}
that multiscale limits with a non-trivial cell structure typically produce only
``simple'' or ``contact'' EDP-convergence, but not tilt-EDP.

This is the regime for which Mielke's theory of non-equilibrium steady states
(NESS) in port gradient systems~\cite{Mielke2023,FrenzelMielke2026} was
developed, and it suggests a second route from the cell formula $\sfM$ to an
$\sfR\oplus\sfR^*$ structure. We sketch it here without striving for rigour.
The central object is the characterization
\begin{equation}\label{eq:M=sup}
    \sfM(\gamma,\jb) = \sup_{\eta\in H^{1/2}(\Gamma)} \set*{ \skp{\eta,\jb}_\Gamma - \sfB(\gamma,\eta)}
\end{equation}
with the \emph{reduced BER functional}
\begin{equation}\label{eqdef:sfB}
    \sfB(\gamma,\eta) = \sup_{u\in \Acc:\,\tr u = \ub}\ \inf_{\xi\in H^1(\Omega,u\pi):\,\tr\xi=\eta}\ \set*{\frac{1}{2} \int_\Omega \abs*{\nabla \xi}^2 u \dx\pi - \calI(u)} \,.
\end{equation}
Indeed, by Lemma~\ref{lem:M:R:Rstar:char} and the dualization of the weighted
negative half-order norm,
\[
\sfM(\gamma,\jb) = \inf_{u: \tr u=\ub} \sup_{\eta\in H^{1/2}(\Gamma)}\;\set*{ \skp*{\jb,\eta}_\Gamma - \frac{1}{2} \norm{\eta}_{\dot H^{1/2}(\Gamma,u \pi)}^2 + \calI(u) } \,,
\]
so that~\eqref{eq:M=sup} becomes equivalent to showing the absence of a duality gap, i.e., that the infimum over $u$ and the supremum over $\eta$ can be exchanged.
We leave this for future work.

The functional $\sfB$ arises naturally in the fast--slow reading of the
$\eps$-problem, in which all quantities on $\Omega$ are fast and all quantities
on $\Gamma$ are slow. The BER functional of the $\eps$-gradient system is
\[
\calB_{\calE^\eps,\calR^\eps}\bra[\big]{(\mu,\gamma),(\Xi,\zeta)}
\coloneqq (\calR^\eps)^*\bra[\big]{(\mu,\gamma),(\Xi,\zeta)}
- (\calR^\eps)^*\bra[\big]{(\mu,\gamma),\,-\anabla\rmD\calE^\eps(\mu,\gamma)}\, ,
\]
and $\calD^\eps_T$ is recovered from it by dualization,
\[
    \calD^\eps_T\bra[\big]{(\mu,\gamma),(J,j)} = \int_0^T
        \sup_{ \Xi,\zeta}\set*{\int_\Omega\Xi\cdot\dx J+\int_\Gamma\zeta\dx j - \calB_{\calE^\eps,\calR^\eps}\bra[\big]{(\mu,\gamma),(\Xi,\zeta)} } \dx t .
\]
Evaluating $\calB_{\calE^\eps,\calR^\eps}$ on the test pairs
$\anabla(\xi,\eta)=(\nabla\xi,\eta-\tr\xi)$ with $\tr\xi=\eta$ makes both
boundary contributions vanish, because on the constraint set we have $\tr u=\ub$ so that the
boundary slot of $-\anabla\rmD\calE^\eps$ is $\log(\tr u/\ub)=0$.
We observe that the two interior contributions reduce 
to the two terms in~\eqref{eqdef:sfB} even for $\eps>0$.
Hence, taking the supremum over $\xi$ and the infimum over $u$ gives $\sfB$ as lower bound in the limit. 
Finally, the $\eps^{-1/2}$ prefactor of the boundary terms 
enforces $\eta=\tr\xi$ and $\ub=\tr u$ in the limit, 
which gives rise to the two constraints
in~\eqref{eqdef:sfB}.
A rigorous version of this argument is also left as future work.

\subsection*{What changes without the nondegeneracy assumption}
\label{sec:discussion:degenerate}

The assumption~\eqref{hyp:P} of uniform bounds on $u$ and $\ub$ was used in four places, and each of them
becomes a separate piece of analysis when it is dropped. 
We list them, both as a
guide to future work and as an indication of where the precise difficulty lies.

\begin{enumerate}[label=(\alph*)]
\item \emph{The class of interior densities.} Without the uniform bounds, the natural class is the cone
$\FI(\Omega,\pi)\coloneqq\set{u\in L^1_{\geq0}(\Omega,\pi):\sqrt u\in H^1(\Omega,\pi)}$
on which $\calI$ is finite, and the trace has to be replaced by the
\emph{Fisher trace} $\tr u \coloneqq (\tr\sqrt u)^2$, which is positively
one-homogeneous and monotone but \emph{not} linear. This was a crucial ingredient in the convexity part of Lemma~\ref{lem:props:cell}.

\item \emph{The weighted half-order spaces.} 
For a weight $u\pi$ vanishing on part of~$\Gamma$,
a trace of elements of $H^1(\Omega,u\pi)$ can only exist
on the set $\set{\ub>0}$.
Hence, the trace space in~\eqref{eq:def:H12} must be replaced
by a quotient construction.
That is, one works with the homogeneous space $\dot H^1(\Omega,u\pi)$
modulo constants, with $H^1_\ast(\Omega,u\pi)$ defined as
the closure of the smooth functions supported away from $\set{\ub>0}$,
and sets
$\dot H^{1/2}(\Gamma,u\pi)\coloneqq \dot H^1(\Omega,u\pi)/H^1_\ast(\Omega,u\pi)$.
Elements of this space are classes of extensions, not functions on~$\Gamma$,
and the inequality~\eqref{eq:H12:equivalence} is lost.

\item \emph{Optimality of the harmonic extension.} With the larger class of
admissible perturbations coming from~(b), the vanishing of the cross term in
the proof of Lemma~\ref{lem:logu:optimal} is no longer the definition of the
variational solution, but needs a proof.
In addition,
the identity $\tr\calU(\ub)=\ub$ becomes a genuine boundary-regularity
statement for unbounded data.

\item \emph{The temporal reference measure.} 
Without a two-sided bound the
interior mass and the interior flux need not have $L^1$-densities in time, so
the limit continuity equation must be formulated with a general reference
measure $\sigma\in\calM_{\geq0}([0,T])$ dominating $\gamma(\cdot\times\Gamma)$
and $\abs{J}(\cdot\times\Omega)$.
The limit dissipation is a priori
$\int_0^T\sfM(\gamma_t,\jb_t)\dx\sigma(t)$. 
It has to be shown that the singular part of $\sigma$
carries no cost, and hence that $t\mapsto\gamma_t$ is absolutely continuous.
\end{enumerate}
We believe that none of these points changes the structure of the results and we expect that both
Result~\ref{result:EDPlimit} and Result~\ref{result:CELL:splitting} can be formulated under suitable assumptions, e.g., the
latter under the hypothesis $\calI(\calU(\ub))<\infty$ on the datum $\ub$. 

\subsection*{A possible metric formulation}

The cell formula arises from optimizing the interior configuration
in~\eqref{eq:intro:CS}. Fixing there instead the $\pi$-harmonic extension
$\calU(\ub)$ irrespective of $\dot\gamma$, which is an admissible, but in general not the minimal choice, we obtain the formally non-negative functional
\begin{equation}\label{eq:def:DeGiorgi}
    \begin{split}
	\calL^{\DtN,2}_T(\gamma) = &\calH(\gamma_T|\omega) - \calH(\gamma_0|\omega)\\
    & + \frac{1}{2} \int_0^T \mkern-4mu \bra*{ \norm{\dot \gamma}_{\dot H^{-1/2}(\Gamma,\calU(\ub)\pi)}^2
		\! + \norm{\log \ub}_{\dot H^{1/2}(\Gamma,\calU(\ub)\pi)}^2 } \! \dx{t} \,.
    \end{split}
\end{equation}
The functional $\calL^{\DtN,2}_T$ is induced by a formal Riemannian structure,
and we expect it to describe a metric curve of maximal slope with respect to a
nonlocal transport distance on $\calP(\Gamma)$ induced by the metric tensor
$\norm{\cdot}^2_{\dot H^{-1/2}(\Gamma,\calU(\ub)\pi)}$. We are not aware of a
coupled bulk--boundary diffusion having $\calL^{\DtN,2}_T$ as its EDP-limit,
and leave both questions to future work.

\begin{acknowledgement}
The authors are grateful to Delio Mugnolo for pointing us to the
Dirichlet-to-Neumann operator and the question about the derivation of its
gradient structure, and to Wolfgang Arendt and Tom ter Elst for the background on its semigroup
formulation. In addition, we are grateful to Jasper Hoeksema for pointing out the literature on the classical boundary process.
We thank the anonymous referee for the careful reading of the manuscript and for the valuable suggestions, which improved the presentation. \\
GH is funded by the Deutsche Forschungsgemeinschaft (DFG, German Research Foundation) under Germany's Excellence Strategy – The Berlin Mathematics Research Center MATH+ (EXC-2046/2, project ID: 390685689).\\
The research of AS is partially based upon work from COST Action 24122 mSPACE,
supported by COST (European Cooperation in Science and Technology), www.cost.eu.
\medskip
\item[\hskip\labelsep{\bfseries AI Declaration}]
The authors used Claude (Anthropic) for proofreading the manuscript and for
checking it for consistency and language improvements.
In particular, the measurability of the integrand in~\eqref{eqdef:D0} was
pointed out by Claude and subsequently established by the authors from the lower
semicontinuity of Proposition~\ref{prop:M:lsc}.
The authors take full
responsibility for the content of the manuscript.
\end{acknowledgement}

\bibliographystyle{spmpsci}
\bibliography{MultiscaleDirichletNeumann}

\end{document}